\documentclass[a4paper,leqno,draft]{article}
\usepackage{amssymb,amsmath,amsfonts,amsthm,
mathrsfs}

\renewcommand{\theequation}{\thesection.\arabic{equation}}
\newtheorem{theorem}{Theorem}[section]
\newtheorem{corollary}[theorem]{Corollary}
\newtheorem{lemma}[theorem]{Lemma}
\newtheorem{proposition}[theorem]{Proposition}
\newtheorem{definition}[theorem]{Definition}
\theoremstyle{definition}
\newtheorem{remark}[theorem]{Remark}
\newtheorem{example}[theorem]{Example}

\newcommand{\wt}[1]{\widetilde{#1}}

\newcommand{\Cinf}{\ensuremath{\mathcal{C}^\infty}}
\newcommand{\Cinfc}{\ensuremath{\mathcal{C}^\infty_{\text{c}}}}
\newcommand{\D}{\ensuremath{{\cal D}}}
\renewcommand{\S}{\mathscr{S}}
\newcommand{\E}{\ensuremath{{\cal E}}}

\newcommand{\mb}[1]{\ensuremath{\mathbb{#1}}}
\newcommand{\N}{\mb{N}}

\newcommand{\R}{\mb{R}}
\newcommand{\C}{\mb{C}}

\newcommand{\LL}{\mathcal{L}}

\newcommand{\A}{\ensuremath{{\cal A}}}
\newcommand{\G}{\ensuremath{{\cal G}}}

\newcommand{\Gc}{\ensuremath{{\cal G}_\mathrm{c}}}
\newcommand{\Gcinf}{\ensuremath{{\cal G}^\infty_\mathrm{c}}}
\newcommand{\GS}{\G_{{\, }\atop{\hskip-4pt\scriptstyle\S}}\!}
\newcommand{\EM}{\ensuremath{{\cal E}_{M}}}
\newcommand\EcM{\mathcal{E}_{\mathrm{c},M}}
\newcommand\EcMinf{\mathcal{E}^\infty_{\mathrm{c},M}}
\newcommand\Nc{\mathcal{N}_{\mathrm{c}}}

\newcommand{\ES}{\mathcal{E}_{\S}}
\newcommand{\EMinf}{\ensuremath{{\cal E}^\infty_{M}}}
\newcommand{\ESinf}{\mathcal{E}_{\S}^{\infty}}

\newcommand{\Neg}{\mathcal{N}}
\newcommand{\NS}{\mathcal{N}_{\S}}

\newcommand{\Ginf}{\ensuremath{\G^\infty}}

\newcommand{\GSinf}{\G^\infty_{{\, }\atop{\hskip-3pt\scriptstyle\S}}}

\newcommand{\singsupp}{\mathrm{sing\, supp}}
\newcommand{\supp}{\mathrm{supp}}

\newcommand{\ssc}{\mathrm{sc}}

\newfont{\bigmath}{cmr12 at 13pt}

\newfont{\grecomath}{cmmi12 at 15pt}

\newcommand{\val}{\mathrm{v}} 
\newcommand{\esp}{\mathrm{e}}

\newfont{\bl}{msbm10 scaled \magstep2}

\newcommand{\beq}{\begin{equation}}
\newcommand{\eeq}{\end{equation}}

\newcommand{\notmid}{\mid\kern-0.5em\not\kern0.5em}

\newcommand{\eps}{\varepsilon}

\newcommand{\Om}{\Omega}

\renewcommand{\Re}{\ensuremath{\text{Re}}}
\renewcommand{\Im}{\ensuremath{\text{Im}}}

\newcommand{\M}{\mathcal{M}}
\newcommand{\mC}{\mathcal{C}}
\newcommand{\mF}{\mathcal{F}}
\newcommand{\mB}{\mathcal{B}}
\newcommand{\mP}{\mathcal{P}}

\newcommand{\Lb}{\mathcal{L}_{\rm{b}}}

\newcommand{\dslash}{d\hspace{-0.4em}{ }^-\hspace{-0.2em}}

\begin{document}

\title{{\bf Fundamental solutions in the Colombeau framework: applications to solvability and regularity theory}}

\author{Claudia Garetto\footnote{Supported by FWF (Austria), grants T305-N13 and  Y237-N13.} \\ 
Institut f\"ur Technische Mathematik,\\ Geometrie und Bauinformatik,\\
Universit\"at Innsbruck, Austria\\
\texttt{claudia@mat1.uibk.ac.at}\\
}
\maketitle
\begin{abstract} 
In this article we introduce the notion of fundamental solution in the Colombeau context as an element of the dual $\LL(\Gc(\R^n),\wt{\C})$. After having proved the existence of a fundamental solution for a large class of partial differential operators with constant Colombeau coefficients, we investigate the relationships between fundamental solutions in $\LL(\Gc(\R^n),\wt{\C})$, Colombeau solvability and $\G$- and $\Ginf$-hypoellipticity respectively. 
\end{abstract}

\setcounter{section}{-1}
\section{Introduction}
The purpose of this paper is to address a main question in the theory of partial differential operators with constant Colombeau coefficients: what is a good notion of fundamental solution in this setting?\\ 
With the adjective good we intend to look for a definition able to provide a successful tool of investigation for issues as $\G$- and $\Ginf$-hypoellipticity and Colombeau solvability. As it is common in the recent research work within Colombeau theory (see \cite{Garetto:04, Garetto:05a, Garetto:05b, Garetto:06a, Garetto:06b, GGO:03, GH:05}), we aim to develop a new set-up of concepts and properties by means of which to achieve statements modelled on well-known classical results of distribution theory. More precisely, given a partial differential operator $P(D)=\sum_{|\alpha|\le m}c_\alpha D^\alpha$ with coefficients in the ring $\wt{\C}$ of generalized numbers, we want a notion of fundamental solution $E$ such that the $\G$- and the $\Ginf$-hypoellipticity of $P(D)$ may be understood by looking at $E$ outside the origin and such that a solution $u$ in the Colombeau algebra $\G(\R^n)$ of the equation $P(D)u=v$ may be found as the convolution product $E\ast v$ when $v$ has compact support. 

Differently from some previous attempts due to Pilipovic and alii \cite{NP:98,NPS:98} we settle ourselves in the dual of the Colombeau algebra $\Gc(\R^n)$ instead than in the usual Colombeau algebra $\G(\R^n)$. This means to take the canonical embedding $\iota_d(\delta)$ of the distributional delta into the dual $\LL(\Gc(\R^n),\wt{\C})$ and to consider the equation $P(D)E=\iota_d(\delta)$ in the dual context. It follows that for the first time a fundamental solution $E$ of the operator $P(D)$ is defined not as a generalized function in $\G(\R^n)$ but as a functional in $\LL(\Gc(\R^n),\wt{\C})$. The results achieved in the paper prove that this is a good notion of fundamental solution. 

In the sequel we describe the contents of the sections in more detail. 

Section \ref{section_preliminaries} collects some preliminaries of Colombeau and duality theory. In view of the techniques which will be employed in the sequel, we focus our attention on the convolution product between Colombeau generalized functions and functionals and on the Fourier-Laplace transform in the algebra $\Gc(\R^n)$ of generalized functions with compact support.  

Section \ref{section_fund} is the mathematical core of the paper. Making use of some results due to H\"ormander \cite[Chapter III]{Hoermander:63}, \cite[Chapter X]{Hoermander:V2}, we prove an adapted version of the Malgrange-Ehrenpreiss Theorem for fundamental solutions in the dual $\LL(\Gc(\R^n),\wt{\C})$. More precisely, we prove that under a certain assumption of invertibility in a point of the generalized weight function $\wt{P}(\xi):=\big(\sum_\alpha |\partial^\alpha P(\xi)|^2\big)^{\frac{1}{2}}$, the corresponding operator $P(D)$ admits a fundamental solution $E$ in $\LL(\Gc(\R^n),\wt{\C})$ which can be defined by a moderate net $(E_\eps)_\eps$ of distributions. This kind of net cannot be regarded as a representative of a Colombeau generalized function but becomes meaningful in the dual context as a functional with ``basic structure''. The notion of basic functional has been introduced in \cite[Section 1]{Garetto:06a} and turns out to be crucial in many technical issues concerning regularity theory and microlocal analysis (see \cite{Garetto:06a, Garetto:06b}). As a straightforward application of the previous existence theorem we investigate the solvability of the equation $P(D)u=v$ in $\LL(\Gc(\R^n),\wt{\C})$ when the right-hand side is compactly supported. The family of evolution operators with respect to the halfspace $H_n=\{x\in\R^n:\, x_n\ge 0\}$ is the last topic of Section \ref{section_fund}. In Subsection \ref{subsec_evol} we provide a condition on the generalized polynomial of $P(D)$ which is sufficient to claim that $P(D)$ is an evolution operator with respect to $H_n$ and we discuss some explanatory examples. Our interest for fundamental solutions in $\LL(\Gc(\R^n),\wt{\C})$ supported in a certain halfspace is motivated by the desire of developing in the future a theory of generalized hyperbolic operators (with constant Colombeau coefficients) based on the support's properties of the corresponding fundamental solutions.  

Section \ref{section_hyp} shows that the $\G$- and $\Ginf$-hypoellipticity of a partial differential operator $P(D)$ with coefficients in $\wt{\C}$ may be characterized by making use of the fundamental solutions. In analogy with the classical theory of operators with constant coefficients we obtain that $P(D)$ is $\G$-hypoelliptic if and only if it admits a fundamental solution $E\in\LL(\Gc(\R^n),\wt{\C})$ with basic structure which belongs to $\G$ outside the origin. The same assertion holds by replacing $\G$ with $\Ginf$. After having introduced a notion of $\G$- and $\Ginf$-ellipticity by means of different invertibility conditions on the principal symbol we employ the new fundamental solution methods in proving that ellipticity implies hypoellipticity in our generalized setting.

The recent investigation of the $\G$- and $\Ginf$-regularity properties of generalized differential and pseudodifferential operators in the Colombeau context \cite{Garetto:04, GGO:03, GH:05, GH:05b, HO:03, HOP:05} has provided several sufficient conditions of $\G$- and $\Ginf$-hypoellipticity, i.e. hypotheses on the generalized symbol of the operator $P(D)$ which allow to conclude that a basic functional $T\in\LL(\Gc(\R^n),\wt{\C})$ is actually a generalized function in $\G(\R^n)$ or $\Ginf(\R^n)$ when $P(D)T$ belongs to $\G(\R^n)$ or $\Ginf(\R^n)$ respectively. The search for necessary conditions for $\G$- and $\Ginf$-hypoellipticity has been a long-standing open problem. A necessary condition for $\Ginf$-hypoellipticity on the symbol of a partial differential operator with generalized constant coefficients has been obtained for the first time by the author in \cite{Garetto:06b}, by means of some functional analytic methods involving the closed graph theorem for Fr\'echet $\wt{\C}$-modules. Since these methods cannot be directly applied to the $\G$-hypoellipticity case this part of the necessary conditions' problem has been open so far. In this paper, making use of the characterizations of $\G$- and $\Ginf$-hypoellipticity which come from the existence of a fundamental solution in $\LL(\Gc(\R^n),\wt{\C})$ and of the Fourier-Laplace transform defined on $\Gc(\R^n)$ we achieve a necessary condition for hypoellipticity in both the $\G$- and $\Ginf$-cases. This result involves partial differential operators $P(D)$ whose weight function $\wt{P}$ is invertible in some point of $\R^n$. The necessary condition for $\Ginf$-hypoellipticity obtained in this paper coincides with the one formulated in \cite{Garetto:06b} even though the methods employed in the proofs are completely independent. 

Some interesting examples of fundamental solutions in $\LL(\Gc(\R^n),\wt{\C})$ are collected in Section \ref{section_ex}. Among them we consider a distributional fundamental solution of the operator $(\partial_1...\partial_n)^k$ and we derive a structure theorem for basic functionals in the duals $\LL(\Gc(\R^n),\wt{\C})$ and $\LL(\G(\R^n),\wt{\C})$. For the sake of completeness and the advantage of the reader the paper ends with an appendix on the solvability of the equation $P(D)u=v$ when $v$ is a basic functional in $\LL(\Gc(\R^n),\wt{\C})$. Inspired by the theory of $B_{p,k}$ spaces developed by H\"ormander we investigate the solution $u$ more deeply than in Section \ref{section_fund}, pointing out some specific moderateness properties.

\section{Preliminaries notions}
\label{section_preliminaries}
This section provides some background of Colombeau and duality theory for the techniques employed in the paper. Particular attention is given to the convolution product between Colombeau generalized functions and functionals and to the Fourier-Laplace transform of a Colombeau generalized function with compact support. 

Before dealing with the major points of the Colombeau construction we recall that the regularity issues discussed in Section \ref{section_hyp} will make use of the following concept of \emph{slow scale net (s.s.n)}. A slow scale net is a net $(r_\eps)_\eps\in\C^{(0,1]}$ such that 
\[
\forall q\ge 0\, \exists c_q>0\, \forall\eps\in(0,1]\qquad\qquad\qquad\qquad |r_\eps|^q\le c_q\eps^{-1}.
\]

\subsection{Colombeau generalized functions and duality theory}
As pointed out in \cite{Garetto:05a, Garetto:05b, Garetto:04th, GHO:06} the most common spaces and algebras of generalized functions of Colombeau type can be introduced and investigated under a topological point of view by making use of the following models.

Let $E$ be a locally convex topological vector space topologized through the family of seminorms $\{p_i\}_{i\in I}$. The elements of 
\[
\begin{split} 
\M_E &:= \{(u_\eps)_\eps\in E^{(0,1]}:\, \forall i\in I\,\, \exists N\in\N\quad p_i(u_\eps)=O(\eps^{-N})\, \text{as}\, \eps\to 0\},\\
\M^\ssc_E &:=\{(u_\eps)_\eps\in E^{(0,1]}:\, \forall i\in I\,\, \exists (\omega_\eps)_\eps\, \text{s.s.n.}\quad p_i(u_\eps)=O(\omega_\eps)\, \text{as}\, \eps\to 0\},\\
\M^\infty_E &:=\{(u_\eps)_\eps\in E^{(0,1]}:\, \exists N\in\N\,\, \forall i\in I\quad p_i(u_\eps)=O(\eps^{-N})\, \text{as}\, \eps\to 0\},\\
\Neg_E &:= \{(u_\eps)_\eps\in E^{(0,1]}:\, \forall i\in I\,\, \forall q\in\N\quad p_i(u_\eps)=O(\eps^{q})\, \text{as}\, \eps\to 0\},
\end{split}
\]
 
are called $E$-moderate, $E$-moderate of slow scale type, $E$-regular and $E$-negligible, respectively. We define the space of \emph{generalized functions based on $E$} as the factor space $\G_E := \M_E / \Neg_E$. 

The rings $\wt{\C}=\EM/\Neg$ of \emph{complex generalized numbers} and $\wt{\R}$ of \emph{real generalized numbers} are obtained by taking $E=\C$ and $E=\R$ respectively. $\wt{\R}$ can be endowed with some more structure by defining the order relation: $r\le s$ if and only if there are representatives $(r_\eps)_\eps$, $(s_\eps)_\eps$ with $r_\eps\le s_\eps$ for all $\eps$. It follows that $r\in\wt{\R}$ is positive ($r\ge 0$) if there exists a representative $(r_\eps)_\eps$ such that $r_\eps\ge 0$ for all $\eps\in(0,1]$. An element $r$ of $\wt{\R}$ is called strictly nonzero if there exists some representative $(r_\eps)_\eps$ and an $m\in\N$ such that $|r_\eps|\ge\eps^m$ for all sufficiently small $\eps$. Finally a positive and strictly nonzero $r\in\wt{\R}$ is called strictly positive. This means that $r_\eps\ge \eps^m$ for some representative $(r_\eps)_\eps$, some $m\in\N$ and for all $\eps$ small enough.

For any locally convex topological vector space $E$ the space $\G_E$ has the structure of a $\wt{\C}$-module. The ${\C}$-module $\G^\ssc_E:=\M^\ssc_E/\Neg_E$ of \emph{generalized functions of slow scale type} and the $\wt{\C}$-module $\Ginf_E:=\M^\infty_E/\Neg_E$ of \emph{regular generalized functions} are subrings of $\G_E$ with more refined assumptions of moderateness at the level of representatives. We use the notation $u=[(u_\eps)_\eps]$ for the class $u$ of $(u_\eps)_\eps$ in $\G_E$. This is the usual way adopted in the paper to denote an equivalence class.

The family of seminorms $\{p_i\}_{i\in I}$ on $E$ determines a \emph{locally convex $\wt{\C}$-linear} topology on $\G_E$ (see \cite[Definition 1.6]{Garetto:05a}) by means of the \emph{valuations}
\[
\val_{p_i}([(u_\eps)_\eps]):=\val_{p_i}((u_\eps)_\eps):=\sup\{b\in\R:\qquad p_i(u_\eps)=O(\eps^b)\, \text{as $\eps\to 0$}\}
\] 
and the corresponding \emph{ultra-pseudo-seminorms} $\{\mP_i\}_{i\in I}$. The theoretical presentation concerning  definitions and properties of valuations and ultra-pseudo-seminorms in the abstract context of $\wt{\C}$-modules is here omitted for the sake of brevity and can be found in \cite[Subsections 1.1, 1.2]{Garetto:05a}. 

In the current paper the valuation and the ultra-pseudo-norm on $\wt{\C}$ obtained through the absolute value in $\C$ are denoted by $\val$ and $|\cdot|_{\esp}$ respectively. The Colombeau algebra $\G(\Om)=\EM(\Om)/\Neg(\Om)$ is obtained as a ${\wt{\C}}$-module of $\G_E$-type by choosing $E=\E(\Om)$. The seminorms $p_{K,i}(f)=\sup_{x\in K, |\alpha|\le i}|\partial^\alpha f(x)|$, where $K$ is a compact subset of $\Om$, generate the family of ultra-pseudo-seminorms  $\mP_{K,i}(u)=\esp^{-\val_{p_{K,i}}(u)}$ and give to $\G(\Om)$ the topological structure of a Fr\'echet $\wt{\C}$-module. We recall that $\Om\to\G(\Om)$ is a fine sheaf of differential algebras on $\R^n$ and that the constants of $\G(\R^n)$ are the elements of $\wt{\C}$.

The Colombeau algebra $\Gc(\Om)$ of generalized functions with compact support is topologized by means of a strict inductive limit procedure. More precisely, setting $\G_K(\Om):=\{u\in\Gc(\Om):\, \supp\, u\subseteq K\}$ for $K\Subset\Om$, $\Gc(\Om)$ is the strict inductive limit of the sequence of locally convex topological $\wt{\C}$-modules $(\G_{K_n}(\Om))_{n\in\N}$, where $(K_n)_{n\in\N}$ is an exhausting sequence of compact subsets of $\Om$ such that $K_n\subseteq K_{n+1}$. We recall that the space $\G_K(\Om)$ is endowed with the topology induced by $\G_{\mathcal{D}_{K'}(\Om)}$ where $K'$ is a compact subset containing $K$ in its interior. In detail we consider on $\G_K(\Om)$ the ultra-pseudo-seminorms $\mP_{\G_K(\Om),n}(u)=\esp^{-\val_{K,n}(u)}$. Note that the valuation $\val_{K,n}(u):=\val_{p_{K',n}}(u)$ is independent of the choice of $K'$ when acts on $\G_K(\Om)$. 

Regularity theory in the Colombeau context as initiated in \cite{O:92} is based on the subalgebra $\Ginf(\Om)=\EMinf(\Om)/\Neg(\Om)$ of $\G(\Om)$ obtained as $\Ginf_E$-space when $E=\E(\Om)$. The intersection of $\Ginf(\Om)$ with $\Gc(\Om)$ defines $\Gcinf(\Om)$. We finally consider the Colombeau algebras $\GS(\R^n)=\ES(\R^n)/\NS(\R^n)$ and $\GSinf(\R^n)=\ESinf(\R^n)/\NS(\R^n)$ of generalized functions based on $\S(\R^n)$ determined as $\G_E$ and $\Ginf_E$ spaces respectively by taking $E=\S(\R^n)$. From a topological point of view $\GS(\R^n)$ and $\Ginf(\Om)$ are Fr\'echet $\wt{\C}$-modules, $\Gcinf(\Om)$ is the strict inductive limit of a family of ultra-pseudo-normed $\wt{\C}$-modules and $\GSinf(\R^n)$ is an ultra-pseudo-normed $\wt{\C}$-module.

A duality theory for $\wt{\C}$-modules had been developed in \cite{Garetto:05a, Garetto:05b} in the framework of topological and locally convex topological $\wt{\C}$-modules. Starting from an investigation of $\LL(\G,\wt{\C})$, the $\wt{\C}$-module of all $\wt{\C}$-linear and continuous functionals on $\G$, it provides the theoretical tools for dealing with the topological duals of the Colombeau algebras $\Gc(\Om)$, $\G(\Om)$ and $\GS(\R^n)$. The spaces $\LL(\G(\Om),\wt{\C})$, $\LL(\Gc(\Om),\wt{\C})$ and $\LL(\GS(\R^n),\wt{\C})$ are endowed with the \emph{topology of uniform convergence on bounded subsets} (see \cite[Remark 2.11]{Garetto:05a}) and, as proven in \cite[Theorems 3.1, 3.8]{Garetto:05b}, the following chains of continuous embeddings hold:  
\beq
\label{chain_1}
\Ginf(\Om)\subseteq\G(\Om)\subseteq\LL(\Gc(\Om),\wt{\C}),
\eeq
\[
\Gcinf(\Om)\subseteq\Gc(\Om)\subseteq\LL(\G(\Om),\wt{\C}),
\]
\[
\GSinf(\R^n)\subseteq\GS(\R^n)\subseteq\LL(\GS(\R^n),\wt{\C}).
\]
Since $\Om\to\LL(\Gc(\Om),\wt{\C})$ is a sheaf we can define the \emph{support of a functional $T$} (denoted by $\supp\, T$). In analogy with distribution theory from Theorem 1.2 in \cite{Garetto:05b} we have that $\LL(\G(\Om),\wt{\C})$ can be identified with the set of functionals in $\LL(\Gc(\Om),\wt{\C})$ having compact support. 

The Colombeau algebra $\GS(\R^n)$ and its dual $\LL(\GS(\R^n),\wt{\C})$ are the natural setting where to define the Fourier transform $\mF$ and its inverse $\mF^{-1}$. In detail we employ the classical definition at the level of representatives in $\GS(\R^n)$ and the definition $\mF(T)(u)=T(\mF(u))$ on the functionals $T$ of $\LL(\GS(\R^n),\wt{\C})$. The reader may refer to \cite[Subsection 1.4]{Garetto:06a} for further explanation. Since $\Gc(\Om)\subseteq\GS(\R^n)$ we are already able to compute the Fourier transform of a Colombeau generalized function with compact support and we will extend $\mF:\Gc(\Om)\to\GS(\R^n)$ to the Fourier-Laplace transform $\mF\LL$ in Subsection \ref{subsec_FL}.  

As already observed in \cite{Garetto:06a, GHO:06}, the chains of inclusions in \eqref{chain_1} make it meaningful to measure the regularity of a functional in $\LL(\Gc(\Om),\wt{\C})$ with respect to the algebras $\G(\Om)$ and $\Ginf(\Om)$. We define the \emph{$\G$-singular support} of $T$ (${\rm{singsupp}}_\G\, T$) as the complement of the set of all points $x\in\Om$ such that the restriction of $T$ to some open neighborhood $V$ of $x$ belongs to $\G(V)$. 
Analogously replacing $\G$ with $\Ginf$ we introduce the notion of \emph{$\Ginf$-singular support} of $T$ denoted by ${\rm{singsupp}}_{\Ginf} T$. A microlocal analysis in the double $\G$- and $\Ginf$-version has been developed in the dual $\LL(\Gc(\Om),\wt{\C})$ by making use of the notions of $\G$- and $\Ginf$-wave front set \cite{Garetto:06a}. In this context a main role is played by the functionals in $\LL(\Gc(\Om),\wt{\C})$ and $\LL(\G(\Om),\wt{\C})$ which have a ``basic'' structure. In detail, we say that $T\in\LL(\Gc(\Om),\wt{\C})$ is basic if there exists a net $(T_\eps)_\eps\in\D'(\Om)^{(0,1]}$ fulfilling the following condition: for all $K\Subset\Om$ there exist $j\in\N$, $c>0$, $N\in\N$ and $\eta\in(0,1]$ such that
\beq
\label{mod_basic_1}
\forall f\in\Cinf_K(\Om)\, \forall\eps\in(0,\eta]\qquad\quad
|T_\eps(f)|\le c\eps^{-N}\sup_{x\in K,|\alpha|\le j}|\partial^\alpha f(x)|
\eeq
and $Tu=[(T_\eps u_\eps)_\eps]$ for all $u\in\Gc(\Om)$. For shortness we denote the set of nets of distributions fulfilling the property \eqref{mod_basic_1} by $\M(\Cinfc(\Om),\C)$.\\
In the same way a functional $T\in\LL(\G(\Om),\wt{\C})$ is said to be basic if there exists a net  $(T_\eps)_\eps\in\E'(\Om)^{(0,1]}$ such that there exist $K\Subset\Om$, $j\in\N$, $c>0$, $N\in\N$ and $\eta\in(0,1]$ with the property 
\[
\forall f\in\Cinf(\Om)\, \forall\eps\in(0,\eta]\qquad\quad
|T_\eps(f)|\le c\eps^{-N}\sup_{x\in K,|\alpha|\le j}|\partial^\alpha f(x)|
\]
and $Tu=[(T_\eps u_\eps)_\eps]$ for all $u\in\G(\Om)$.\\
Clearly the sets $\Lb(\Gc(\Om),\wt{\C})$ and $\Lb(\G(\Om),\wt{\C})$ of basic functionals are $\wt{\C}$-linear subspaces of $\LL(\Gc(\Om),\wt{\C})$ and $\LL(\G(\Om),\wt{\C})$ respectively.

\subsection{Convolution between generalized functions and functionals}
We recall some of the properties of the convolution product between functionals and Colombeau generalized functions which are employed in the course of the paper. We refer for definitions and proofs to \cite{Garetto:06a}.
\begin{proposition}
\label{prop_conv_1}
The $\wt{\C}$-bilinear map
\[
(S,T)\to S\ast T:=S_xT_y(\cdot(x+y))
\]
\begin{itemize}
\item[(i)] from $\Gc(\R^n)\times\Lb(\Gc(\R^n),\wt{\C})$ into $\G(\R^n)$,
\item[(ii)] from $\G(\R^n)\times\Lb(\G(\R^n),\wt{\C})$ into $\G(\R^n)$,
\item[(iii)] from $\GS(\R^n)\times\Lb(\G(\R^n),\wt{\C})$ into $\GS(\R^n)$,
\item[(iv)] from $\Gcinf(\R^n)\times\Lb(\Gc(\R^n),\wt{\C})$ into $\Ginf(\R^n)$,
\item[(v)] from $\Ginf(\R^n)\times\Lb(\G(\R^n),\wt{\C})$ into $\Ginf(\R^n)$,
\item[(vi)] from $\GSinf(\R^n)\times\Lb(\G(\R^n),\wt{\C})$ into $\GSinf(\R^n)$,
\item[(vii)] from $\LL(\G(\R^n),\wt{\C})\times\Lb(\Gc(\R^n),\wt{\C})$ into $\LL(\Gc(\R^n),\wt{\C})$,
\item[(viii)] from $\LL(\Gc(\R^n),\wt{\C})\times\Lb(\G(\R^n),\wt{\C})$ into $\LL(\Gc(\R^n),\wt{\C})$,
\item[(ix)] from $\LL(\G(\R^n),\wt{\C})\times\Lb(\G(\R^n),\wt{\C})$ into $\LL(\G(\R^n),\wt{\C})$
\end{itemize}
is separately continuous. Moreover, when at least one of the functionals $S$ and $T$ has compact support and $T$ is basic, the inclusion 
\[
\supp(S\ast T)\subseteq \supp\, S+\supp\, T
\]
holds. If also the functional $S$ is basic then 
\[
\singsupp_\G(S\ast T)\subseteq \singsupp_\G\, S+\singsupp_\G\, T
\]
and
\[
\singsupp_{\Ginf}(S\ast T)\subseteq \singsupp_{\Ginf}\, S+\singsupp_{\Ginf}\, T.
\]
\end{proposition}
\begin{proof}
We only prove the two final inclusions concerning the $\G$- and the $\Ginf$-singular supports. Assume that $T$ has  compact support and take $\psi\in\Cinf(\R^n)$ identically $1$ in a neighborhood of $\singsupp_\G\, T$. Then, we can write $T=T_1+T_2$ with $T_1:=\psi T\in\Lb(\G(\R^n),\wt{\C})$ and $T_2:=(1-\psi)T\in\Gc(\R^n)$. It follows that $S\ast T_2\in\G(\R^n)$ while $S\ast T_1$ is a generalized function on the open set $\{x:\, x-\supp\,T_1\subseteq\R^n\setminus\singsupp_\G\, S\}$. This means that 
\[
\singsupp_\G(S\ast T)=\singsupp_\G(S\ast T_1)\subseteq \singsupp_\G\, S+\supp T_1\subseteq \singsupp_\G\, S+\supp\, \psi.
\]
Since $\supp\, \psi$ can be taken as close to $\singsupp_\G\, T$ as we wish, we obtain the desired inclusion. The proof of the assertion with the $\Ginf$-singular supports is analogous and left to the reader.
\end{proof}
Note that the convolution of two basic functionals is a basic functional too. 
The convolution $S\ast T$ can be defined in many cases when neither $S$ nor $T$ has compact support. What we need is the proper map condition.
\begin{proposition}
\label{prop_conv_proper}
Let $S\in\LL(\Gc(\R^n),\wt{\C})$ and $T\in\Lb(\Gc(\R^n),\wt{\C})$ such that the map
\[
\mu:\supp\,S\times\supp\,T\to \R^n:(x,y)\to x+y
\]
is proper. Then the convolution $S\ast T$ can be defined as a functional in $\LL(\Gc(\R^n),\wt{\C})$. Furthermore, $\supp(S\ast T)\subseteq \supp\,S+\supp\,T$.  
\end{proposition}
\begin{proof}
Let $(V_k)_k$ be an open covering of $\R^n$ such that $V_{k-1}\subseteq\overline{V_k}\Subset V_{k+1}$. By the hypothesis on $\mu$ we know that $\pi_1(\mu^{-1}(\overline{V_k}))$ and $\pi_2(\mu^{-1}(\overline{V_k}))$ are compact subsets of $\R^n$. Let $\phi_{k,1}, \phi_{k,2}\in\Cinfc(\R^n)$ identically $1$ in a neighborhood of $\pi_1(\mu^{-1}(\overline{V_k}))$ and $\pi_2(\mu^{-1}(\overline{V_k}))$ respectively. Then we can define the convolution product $\phi_{k,1}S\ast \phi_{k,2}T|_{V_k}$ as an element of $\LL(\Gc(V_k),\wt{\C})$. This functional does not depend on the choice of the cut-off functions $\phi_{k,1}$ and $\phi_{k,2}$. Indeed, given $\phi_{k,1}$ and $\phi_{k,2}$ with the same neighborhood-property we can write $\phi_{k,1}S\ast \phi_{k,2}T|_{V_k}-\phi'_{k,1}S\ast \phi'_{k,2}T|_{V_k}$ as 
\[
(\phi_{k,1}-\phi'_{k,1})S\ast\phi_{k,2}T|_{V_k}+\phi'_{k,1}S\ast(\phi_{k,2}-\phi'_{k,2})T|_{V_k},
\] 
where both the summand are null. So, we can set $(S\ast T)_k:=\phi_{k,1}S\ast \phi_{k,2}T|_{V_k}$ and since $\LL(\Gc(\R^n),\wt{\C})$ is a sheaf it is enough to prove that the family $\{(S\ast T)_k\}_{k\in\N}$ is coherent in order to conclude that it uniquely defines a functional $S\ast T$ in $\LL(\Gc(\R^n),\wt{\C})$. One easily sees that if $k<k'$ then $(S\ast T)_{k'}|_{V_k\cap V_{k'}}=(S\ast T)_{k'}|_{V_k}=\phi'_{k,1}S\ast\phi'_{k,2}T|_{V_k}=\phi_{k,1}S\ast\phi_{k,2}T|_{V_k}=(S\ast T)_k|_{V_k\cap V_{k'}}$. The inclusion $\supp(S\ast T)\subseteq\supp\,S+\supp\,T$ easily follows from the analogous inclusion in Proposition  	\ref{prop_conv_1} and the definition of $S\ast T$.
\end{proof}
The following corollary is obtained by combining the previous definition of product of convolution with Proposition \ref{prop_conv_1}$(ii)$.
\begin{corollary}
\label{gen_func_case}
Let $S\in\LL(\Gc(\R^n),\wt{\C}$ and $u\in\G(\R^n)$ such that the map
\[
\mu:\supp\,S\times\supp\,u\to \R^n:(x,y)\to x+y
\]
is proper. Then, $S\ast u\in\G(\R^n)$.
\end{corollary}
\begin{remark}
\label{rem_cone_case}
Let $\Gamma\subseteq\R^n$ be a closed convex cone which is \emph{proper} in the sense that it does not contain any straight line. The corresponding map $\mu:\Gamma\times\Gamma\to\R^n:(x,y)\to x+y$ is proper (as it is proved in \cite{Hoermander:V1}, p.104). Hence, the convolution makes the set $\{T\in\Lb(\Gc(\R^n),\wt{\C}):\, \supp\, T\subseteq \Gamma\}$ an algebra. Finally, assume that $\Gamma$ is a closed cone contained in $H_n:=\{x\in\R^n:\, x_n\ge 0\}$ such that $\Gamma\cap\{x:\, x_n=0\}=\{0\}$. Then, the map $\mu:\Gamma\times H_n\to\R^n:(x,y)\to x+y$ is proper. Indeed, given the bounded set $\{(x,y):\, |x+y|\le C\}$ if we suppose that there exist sequences $(x_n)$, $(y_n)$ such that $|x_n+y_n|\le C$ with $|x_n|\to \infty$ passing to subsequences we get that $x_n/|x_n|\to x$ and $y_n/|y_n|\to -x$ for some $x\in\Gamma$ and some $-x\in H_n$. Hence, $x=0$ which contradicts $x_n/|x_n|\to x$.  
\end{remark} 
We finally consider the action of a partial differential operator with constant Colombeau coefficients on the convolution of two functionals.
\begin{proposition}
\label{conv_1}
If $P(D)$ is a partial differential operators with coefficients in $\wt{\C}$, $S\in\LL(\G(\R^n),\wt{\C})$ and $T\in\Lb(\Gc(\R^n),\wt{\C})$ then 
\[
P(D)(S\ast T)=P(D)S\ast T= S\ast P(D)T.
\]
The same equalities hold for $S,T\in\LL(\Gc(\R^n),\wt{\C})$ as in Proposition \ref{prop_conv_proper}.
\end{proposition}
\begin{proof}
If $S$ has compact support then an inspection at the level of representatives shows that the following equalities hold for all $u\in\Gc(\R^n)$:
\begin{multline*}
P(D)S\ast T(u)=(P(D)S)_x(T_y(u(x+y)))=S_x({\,}^tP(D)(T_yu(x+y)))= S_x((P(D)T)_y(u)(x+y))\\
=S\ast P(D)T(u)=S_x(T_y({\,}^tP(D)u(x+y)))=(S\ast T)({\,}^tP(D)u)=P(D)(S\ast T)(u).
\end{multline*} 
We leave to the reader to check that the same result holds for $S\in\LL(\Gc(\R^n),\wt{\C})$ and $T\in\Lb(\Gc(\R^n),\wt{\C})$ satisfying the proper map assumption of Proposition \ref{prop_conv_proper}.
\end{proof}

\subsection{The Fourier-Laplace transform in the Colombeau fra\-me\-work}
\label{subsec_FL}
The purpose of this subsection is to deal with the Fourier-Laplace transform in the Colombeau framework of generalized functions with compact support. The collected material is a partial elaboration of \cite{NP:95, NPS:98, Soraggi:96}.

Let $\A(\C^n)$ be the space of all analytic functions on $\C^n$ and $a\ge 0$. We denote the sets of all nets $(u_\eps)_\eps\in\A(\C^n)^{(0,1]}$ such that 
\[
\forall M\in\N\, \exists N\in\N\, \exists c>0\, \exists\eta\in(0,1]\, \forall\eps\in(0,\eta],\ \forall\zeta\in\C^n\qquad
|u_\eps(\zeta)|\le c\,\eps^{-N}(1+|\zeta|)^{-M}\esp^{a|\rm{Im}(\zeta)|}, 
\]
\[
\exists N\in\N\, \forall M\in\N\, \exists c>0\, \exists\eta\in(0,1]\, \forall\eps\in(0,\eta],\ \forall\zeta\in\C^n\qquad
|u_\eps(\zeta)|\le c\,\eps^{-N}(1+|\zeta|)^{-M}\esp^{a|\rm{Im}(\zeta)|} 
\]
and
\beq
\label{neg_fl_a}
\forall M\in\N\, \forall q\in\N\, \exists c>0\, \exists\eta\in(0,1]\, \forall\eps\in(0,\eta],\ \forall\zeta\in\C^n\qquad
|u_\eps(\zeta)|\le c\,\eps^{q}(1+|\zeta|)^{-M}\esp^{a|\rm{Im}(\zeta)|}
\eeq
by $\E_{\mF\LL,a}(\C^n)$, $\E^\infty_{\mF\LL,a}(\C^n)$ and $\Neg_{\mF\LL,a}(\C^n)$ respectively.
 
We set 
\[
\G_{\mF\LL,a}(\C^n):=\frac{\E_{\mF\LL,a}(\C^n)}{\Neg_{\mF\LL,a}(\C^n)}
\]
and
\[
\Ginf_{\mF\LL,a}(\C^n):=\frac{\E^\infty_{\mF\LL,a}(\C^n)}{\Neg_{\mF\LL,a}(\C^n)}.
\]
The following proposition shows that the classical Fourier-Laplace transform at the level of representatives allows to define the Fourier-Laplace transform of a generalized function in $\Gc(\R^n)$ as an element of some factor space $\G_{\mF\LL,a}(\C^n)$.
\begin{proposition}
\label{prop_fl_gc}
For all $M\in\N$ and $a>0$ there exists a constant $c_{M,a}>0$ such that the inequality 
\beq
\label{est_mod_repr}
|\mF\LL(u)(\zeta)|\le c_{M,a}\,(1+|\zeta|)^{-M}\sup_{|\alpha|\le M, |x|\le a}|\partial^\alpha u(x)|\,\sup_{|x|\le a}\esp^{x\rm{Im}(\zeta)}
\eeq
holds for all $\zeta\in\C^n$ and for all $u\in\Cinf(\R^n)$ with $\supp\, u\subseteq\{x:\, |x|\le a\}$.
\end{proposition}

\begin{definition}
\label{FL_transform}
Let $u\in\Gc(\R^n)$ with $\supp\, u\subseteq\{x:\, |x|<a\}$. The \emph{Fourier-Laplace} transform of $u\in\Gc(\R^n)$ is the generalized function
\[
\mF\LL(u)(\zeta)=\int_{\R^n}\esp^{-ix\zeta}u(x)\,dx
\]
in $\G_{\mF\LL,a}(\C^n)$ obtained by applying the corresponding classical transformation on the representatives of $u$ having support contained in a compact subset $K\subseteq\{x:\, |x|<a\}$ uniformly with respect to the parameter $\eps$.
\end{definition}
The well-definedness of $\mF\LL(u)$ in $\G_{\mF\LL,a}(\C^n)$ is guaranteed by Proposition \ref{prop_fl_gc}.\\
Indeed, if $u\in\Gc(\R^n)$ has support contained in $\{x:\, |x|<a\}$ then it has a representative $(u_\eps)_\eps$ such that $\supp\, u_\eps\subseteq K\Subset\{x:\, |x|<a\}$ for some compact set $K$ and for all $\eps\in(0,1]$. By Proposition \ref{prop_fl_gc} we have that $(\mF\LL(u_\eps))_\eps\in\E_{\mF\LL,a}(\C^n)$. Moreover, when $(u'_\eps)_\eps$ is another representative of $u$ with $\supp\, u_\eps$ contained in a certain compact subset $K'$ of $\{x:\, |x|<a\}$ for all $\eps$, then the difference $(u_\eps-u'_\eps)_\eps$ satisfies \eqref{neg_fl_a} with $a>0$ as above. It follows that $(u_\eps-u'_\eps)_\eps\in\Neg_{\mF\LL,a}(\C^n)$ and that $\mF\LL(u)\in\G_{\mF\LL,a}(\C^n)$.

The following theorem provides a deeper investigation of the properties of $\mF\LL(u)$ and the expected Paley-Wiener type results. For technical reasons we will make use of the subset $\underline{\G}_{\,\mF\LL,a}(\C^n)$ of $\G_{\mF\LL,a}(\C^n)$ obtained by assuming that the estimates which characterize the representatives hold in the whole interval $(0,1]$. In the sequel, $\EcM(\R^n)$ denotes the set of all nets $(u_\eps)_\eps\in\EM(\R^n)$ of smooth functions having support contained in a compact subset of $\R^n$ uniformly with respect to the parameter $\eps$. Analogously one defines $\EcMinf(\R^n)$ and $\Nc(\R^n)$.
\begin{theorem}
\label{theo_Paley_function}
\leavevmode
\begin{itemize}
\item[(i)] If $u\in\Gcinf(\R^n)$ and $\supp\, u\subseteq\{x:\, |x|<a\}$ then $\mF\LL(u)\in\Ginf_{\mF\LL,a}(\C^n)$.
\item[(ii)] If $v\in\underline{\G}_{\,\mF\LL,a}(\C^n)$ then there exists $u\in\Gc(\R^n)$ with $\supp\, u\subseteq\{x:\, |x|\le a\}$ such that
\beq
\label{inv_fl}
(\mF\LL(u_\eps)-v_\eps)_\eps\in\Neg_{\mF\LL,a'}(\C^n),
\eeq
for all representative $(u_\eps)_\eps$ of $u$ in $\EcM(\R^n)$, for all representative $(v_\eps)_\eps$ of $v$ and for all $a'>a$.   
\item[(iii)] If $v\in\underline{\G}^\infty_{\,\mF\LL,a}(\C^n)$ then $(ii)$ holds with $u\in\Gcinf(\R^n)$. 
\end{itemize}
\end{theorem}
\begin{proof}
\leavevmode
$(i)$ Assume that $u\in\Gcinf(\R^n)$. Taking a representative $(u_\eps)_\eps$ of $u$ from \eqref{est_mod_repr} we get that $\mF\LL(u_\eps)_\eps\in\E^\infty_{\mF\LL,a}(\C^n)$. This means that $\mF\LL(u)\in\Ginf_{\mF\LL,a}(\C^n)$. 

$(ii)$ Let $(v_\eps)_\eps$ be a representative of $v$. It has the property
\beq
\label{repr_v}
\forall M\in\N\, \exists N\in\N\, \exists c>0\, \forall\eps\in(0,1]\, \forall\zeta\in\C^n\qquad |v_\eps(\zeta)|\le c\eps^{-N}(1+|\zeta|)^{-M}\esp^{a|{\rm{Im}}\zeta|}.
\eeq
Replacing $\zeta$ with $\xi\in\R^n$ and taking $M=n+1$ in \eqref{repr_v}, we easily see that $(v_\eps)_\eps\in\M_{L^1(\R^n)}$ and therefore 
\beq
\label{u_eps}
u_\eps(x)=(2\pi)^{-n}\int_{\R^n}\esp^{ix\xi}v_\eps(\xi)\, d\xi
\eeq
gives a moderate net of continuous functions on $\R^n$. Analogously the choice of $M=n+1+|\alpha|$ in \eqref{repr_v} for any $\alpha\in\N^n$ makes us conclude that $(u_\eps)_\eps$ is a net of smooth functions and more precisely that $(u_\eps)_\eps\in\EM(\R^n)$. Let $u$ be the generalized function in $\G(\R^n)$ generated by $(u_\eps)_\eps$. We want to prove that $u|_{\{x:\, |x|>a\}}=0$ and that \eqref{inv_fl} holds. Since $v_\eps$ is analytic on $\C^n$ the Cauchy's theorem applied to each variable $\zeta_1,...,\zeta_n$ allows to shift the integration in \eqref{u_eps} into the complex domain and to write
\[
u_\eps(x)=(2\pi)^{-n}\int_{{\rm{Im}}\zeta=\xi_0}\hskip-10pt \esp^{ix\zeta}v_\eps(\zeta)\, d\zeta,
\]  
where $\xi_0$ may be any point in $\R^n$. From \eqref{repr_v} we have that for some constant $C'>0$ (independent of $\xi_0$) the estimate
\beq
\label{est_u_xi_0}
|u_\eps(x)|\le C\eps^{-N}\esp^{a|\xi_0|-x\xi_0}\int_{\R^n}(1+|\xi|)^{-n-1}\,d\xi\le C'\eps^{-N}\esp^{a|\xi_0|-x\xi_0}
\eeq
is valid for all $x\in\R^n$ and for all $\eps\in(0,1]$. Assume that $|x|>a$ and choose $\xi_0=tx/|x|$ with $t>0$ in \eqref{est_u_xi_0}. This yields
\[
|u_\eps(x)|\le C'\eps^{-N}\esp^{t(a-|x|)}.
\]
Letting $t\to +\infty$ we conclude that $u_\eps(x)=0$ for all $\eps$ when $|x|>a$. Hence, $\supp\, u\subseteq\{x:\, |x|\le a\}$. Finally, by construction $\widehat{u_\eps}(\xi)=v_\eps(\xi)$. Since $\widehat{u_\eps}(\xi)$ extends to an analytic function on $\C^n$, by the uniqueness of the analytic continuation one has that $\mF\LL(u_\eps)=v_\eps$ on $\C^n$. It follows that the property \eqref{inv_fl} holds for the generalized function $u$ constructed in this way.

$(iii)$ It is immediate to check that when $v\in\underline{\G}^\infty_{\,\mF\LL,a}(\C^n)$ then the generalized function $u$ with representative $(u_\eps)_\eps$ defined in \eqref{u_eps} by means of $(v_\eps)_\eps\in\E^\infty_{\mF\LL,a}(\C^n)$ belongs to $\Gcinf(\R^n)$.
\end{proof}

\begin{remark}
\label{rem_Paley}
The generalized function $u\in\Gc(\R^n)$ with support contained in $\{x:\, |x|\le a\}$ and such that the second assertion of Theorem \ref{theo_Paley_function} is fulfilled is unique. Indeed, assume that there exists another $u'\in\Gc(\R^n)$ having the same properties. Then, the fact that $(\mF\LL(u_\eps)-v_\eps)_\eps\in\Neg_{\mF\LL,a'}(\C^n)$ and $(\mF\LL(u'_\eps)-v_\eps)_\eps\in\Neg_{\mF\LL,a'}(\C^n)$ leads to $(\mF\LL(u_\eps-u'_\eps))_\eps\in\Neg_{\mF\LL,a'}(\C^n)$. In particular this means that for all $M,q\in\N$ there exists $\eta\in(0,1]$ such that 
\beq
\label{two_repr}
|\widehat{u_\eps}(\xi)-\widehat{u'_\eps}(\xi)|\le \eps^q(1+|\xi|)^{-M},
\eeq
for all $\xi\in\R^n$ and $\eps\in(0,\eta]$. Since, $(\widehat{u_\eps}-\widehat{u'_\eps})_\eps\in\ES(\R^n)$, the characterization of the ideal $\NS(\R^n)$ (see \cite[Proposition 3.4]{Garetto:05b}) allows to conclude from \eqref{two_repr} that $(\widehat{u_\eps}-\widehat{u'_\eps})_\eps\in\NS(\R^n)$. By the injectivity of the Fourier transform on $\GS(\R^n)$ we have that $u=u'$ in $\GS(\R^n)$ and since $\Gc(\R^n)$ is embedded in $\GS(\R^n)$, $u=u'$ in $\Gc(\R^n)$.
\end{remark}
Let now $P(D)$ be a partial differential operator with coefficients in $\wt{\C}$. We leave to the reader to check that the equality 
\beq
\label{eq_FL}
\mF\LL(Pu)(\zeta)=\mF\LL(u)(\zeta)P(\zeta)
\eeq
holds in $\G_{\mF\LL,a}(\C^n)$ for all the partial differential operators $P$ with coefficients in $\wt{\C}$ and for all generalized functions $u\in\Gc(\R^n)$ with $\supp\, u\subseteq\{x:\, |x|<a\}$. Note that in the right hand-side of \eqref{eq_FL} the product between elements of $\G_{\mF\LL,a}(\C^n)$ and polynomials with generalized constant coefficients is defined componentwise at the level of representatives and gives a generalized function in $\G_{\mF\LL,a}(\C^n)$. Since the Fourier-Laplace transform extends the Fourier transform from real variables to complex variables in the classical distributional context as well as in the generalized Colombeau framework, we give to the three notations  $\mF\LL(u)(\zeta)$, $\mF(u)(\zeta)$, $\widehat{u}(\zeta)$ the same meaning. 

We conclude this section of preliminary notions by considering the Fourier-Laplace transform of $u\in\Gc(\R^n)$ computed in the points of $\C^n$ of the form $(\xi',\xi_n+i\eta_n)$, where $\xi'\in\R^{n-1}$ and $\xi_n,\eta_n\in\R$. More in general we will assume that $\eta_n$ is a generalized number in $\wt{\R}$ satisfying suitable logarithmic growth conditions as follows.
\begin{definition}
\label{def_esp_point}
We say that $\xi\in\wt{\R^n}$ is of log-type if there exists a representative $(\xi_\eps)_\eps$ of $\xi\in\wt{\R}^n$ such that $|\xi_\eps|=O(\log(1/\eps))$.
\end{definition} 
The previous condition can equivalently stated saying that the net $(\esp^{|\xi_\eps|})_\eps$ is moderate. As proved in \cite[Subsection 5.1.1]{Garetto:06b} if $\xi\in\wt{\R^n}$ is of log-type then $\esp^{|\xi|}:=(\esp^{|\xi_\eps|})_\eps+\Neg$ and $\esp^{x\xi}:=(\esp^{x\xi_\eps})_\eps+\Neg(\R^n)$ are well-defined elements of $\wt{\R}$ and $\Ginf(\R^n)$ respectively.
 
We are now ready to state the following proposition where we write any $y\in\R^n$ as $(y',y_n)$ with $y'\in\R^{n-1}$ and $y_n\in\R$.  
\begin{proposition}
\label{prop_FL}
\leavevmode
\begin{itemize}
\item[(i)]If $u\in\Gc(\R^n)$ then for all $\eta_n\in\wt{\R}$ of log-type the generalized function 
\beq
\label{def_func}
(\xi',\xi_n)\to\mF\LL(u)(\xi',\xi_n+i\eta_n):=\biggl[\biggl(\int_{\R^n}\esp^{-i(x'\xi'+x_n\xi_n)}\esp^{x_n\eta_{n,\eps}}u_\eps(x',x_n)\, dx'\, dx_n\biggr)_\eps\biggr]
\eeq
belongs to $\GS(\R^n)$. 
\item[(ii)] If $u\in\Gcinf(\R^n)$ then \eqref{def_func} defines a generalized function in $\GSinf(\R^n)$.
\end{itemize}
\end{proposition}
\begin{proof}
Let $u\in\Gc(\R^n)$ with $\supp\, u\subseteq\{x:\, |x|<a\}$. Proposition \ref{prop_fl_gc} proves that the estimate  
\beq
\label{S_estimate}
\biggl\vert\int_{\R^n}\esp^{-ix\zeta}u_\eps(x)\, dx\biggr\vert\le c_{M,a}(1+|\zeta|)^{-M}\sup_{|x|\le a,|\alpha|\le M}|\partial^\alpha u_\eps(x)|\ \esp^{a|\rm{Im}(\zeta)|} 
\eeq
holds for all $(u_\eps)_\eps\in\EcM(\R^n)$ with $\supp\, u_\eps\subseteq\{x:\, |x|<a\}$ and for all $\eps\in(0,1]$. Hence, taking $\zeta=(\xi',\xi_n+i\eta_{n,\eps})$ in \eqref{S_estimate} and computing the derivatives in $(\xi',\xi_n)$ of the net in \eqref{def_func} it follows that 
\begin{multline}
\label{est_deriv}
\biggl\vert\partial^\beta_{\xi'}\frac{d^h}{d\xi_n^h}\int_{\R^n}\esp^{-i(x'\xi'+x_n\xi_n)}\esp^{x_n\eta_{n,\eps}}u_\eps(x',x_n)\, dx'\, dx_n\biggr\vert\le c_{M,a}\,(1+|\xi'|+|\xi_n|)^{-M}\\
\sup_{|(x',x_n)|\le a,|\alpha|\le M}|\partial^\alpha( {x'}^\beta x_n^hu_\eps(x',x_n))|\ \esp^{a|\eta_{n,\eps}|} 
\end{multline}
By the log-type assumption on $\eta_n\in\wt{\R}$ it follows that the net in \eqref{est_deriv} belongs to $\ES(\R^n)$, $\ESinf(\R^n)$ and $\NS(\R^n)$ when $(u_\eps)_\eps$ is an element of $\EcM(\R^n)$, $\EcMinf(\R^n)$ and $\Nc(\R^n)$ respectively. 

In order to conclude that \eqref{def_func} gives a well-defined generalized function in $\GS(\R^n)$ it remains to prove that its definition does not depend on the representative of log-type of $\eta_n$. Let $(\eta_\eps)_\eps$ and $(\eta'_\eps)_\eps$ be representatives of $\eta$ such that $(\esp^{|\eta_\eps|})_\eps)$ and $(\esp^{|\eta'_\eps|})_\eps)$ are moderate. Since we already know that $(\esp^{x_n\eta_{n,\eps}})_\eps\in\EMinf(\R)$ and $(\esp^{x_n\eta_{n,\eps}}-\esp^{x_n\eta'_{n,\eps}})_\eps\in\Neg(\R)$ we can assert that for all $q\in\N$ there exists $\eps_q\in(0,1]$ such that the estimate 
\[
\biggl\vert\int_{\R^n}\esp^{-i(x'\xi'+x_n\xi_n)}\esp^{x_n\eta_{n,\eps}}\biggl(\esp^{x_n\eta_{n,\eps}}-\esp^{x_n\eta'_{n,\eps}}\biggr)u_\eps(x',x_n)\, dx'\, dx_n\biggr\vert\le c\eps^q\sup_{|(x',x_n)|\le a}|u_\eps(x',x_n))|
\]
holds for all $\eps\in(0,\eps^q]$. By the characterization of the ideal $\NS(\R^n)$ it follows that the net given by the integral above belongs to $\NS(\R^n)$. As a consequence \eqref{def_func} defines a generalized function in $\GS(\R^n)$ when $u\in\Gc(\R^n)$ and a generalized function in $\GSinf(\R^n)$ when $u\in\Gcinf(\R^n)$.
\end{proof}
\begin{remark}
\label{rem_FL}
Note that if $\supp\, u_\eps\subseteq\{x\in\R^n:\, x_n>0\}$ for all $\eps\in(0,1]$ then 
\[
\mF\LL(u_\eps)(\xi',\xi_n+i\eta_{n,\eps})=\int_{\R^n}\esp^{-i(x'\xi'+x_n\xi_n)}\esp^{x_n\eta_{n,\eps}}u_\eps(x',x_n)\, dx'\, dx_n
\]
fulfills the following estimate
\[
|\mF\LL(u_\eps)(\xi',\xi_n+i\eta_{n,\eps})|\le c_{M,a}(1+|\xi'|+|\xi_n|)^{-M}\sup_{|(x',x_n)|\le a,|\alpha|\le M}|\partial^\alpha u_\eps(x',x_n)|
\]
for all $(\xi',\xi_n)\in\R^n$, for all $\eta_{n,\eps}<0$ and for all $\eps\in(0,1]$.
\end{remark}
\section{Fundamental solutions in the dual of a Colombeau algebra}
\label{section_fund}
In the sequel $\iota_d$ denotes the embedding of $\D'(\Om)$ into the dual $\LL(\Gc(\Om),\wt{\C})$ given by $\iota_d(w)(u)=[(w(u_\eps))_\eps]$, where $w\in\D'(\Om)$ and $u\in\Gc(\Om)$.
\subsection{Definition}
\begin{definition}
\label{def_fund_sol}
Let $P(D)$ be a partial differential operator with constant Colombeau coefficients. We say that $E\in\LL(\Gc(\R^n),\wt{\C})$ is a fundamental solution of $P(D)$ if $P(D)E=\iota_d(\delta)$ in $\LL(\Gc(\R^n),\wt{\C})$.
\end{definition}
 
\subsection{The Malgrange-Ehrenpreiss Theorem for fundamental solutions in the space $\Lb(\Gc(\R^n),\wt{\C})$}
In this subsection we prove the existence of a fundamental solution in the dual $\LL(\Gc(\R^n),\wt{\C})$ for a large class of partial differential operators with coefficients in $\wt{\C}$. We begin by fixing some language and notations. 
 
Let $P$ be a partial differential operator of order $m$ with coefficients in $\wt{\C}$. Any net of polynomials $(P_\eps)_\eps$ determinated by a choice of representatives of the coefficients of $P$ is called a  representative of $P$. Consider the function $\wt{P}:\R^n\to\wt{\R}$ defined by
\[
\wt{P}^2(\xi)=\sum_{|\alpha|\le m}|\partial^\alpha P(\xi)|^2.
\]
The arguments in \cite[(2.1.10)]{Hoermander:63} yield the following assertion: there exists $C>0$ depending only on $m$ and $n$ such that for all $(P_\eps)_\eps$ the inequality
\beq
\label{est_Hoer}
\wt{P_\eps}(\xi+\eta)\le (1+C|\xi|)^m\wt{P_\eps}(\eta)
\eeq
is valid for all $\xi,\eta\in\R^n$ and all $\eps\in(0,1]$. When the function $\wt{P}:\R^n\to\wt{\R}$ is invertible in some point $\xi_0$ of $\R^n$ Lemma 7.5 in \cite{HO:03} proves that for all representative $(P_\eps)_\eps$ of $P$ there exist $N\in\N$ and $\eta\in(0,1]$ such that
\beq
\label{est_inv}
\wt{P_\eps}(\xi)\ge \eps^N(1+C|\xi_0-\xi|)^{-m},
\eeq
for all $\xi\in\R^n$ and $\eps\in(0,\eta]$. Note that the constant $C>0$ is the same appearing in \eqref{est_Hoer} and $\eps^N$ comes from the invertibility in $\wt{\R}$ of $\wt{P}(\xi_0)$.

We finally recall that $\mathcal{K}$ is the set of tempered weight functions introduced by H\"ormander in
\cite[Definition 2.1.1]{Hoermander:63}, i.e., the set of all positive functions $k$ on $\R^n$ such that for some
constants $C>0$ and $N\in\N$ the inequality 
\[
\label{ineq_weight}
k(\xi+\eta)\le (1+C|\xi|)^N k(\eta) 
\]
holds for all $\xi,\eta\in\R^n$. 
\begin{definition}
\label{def_cla_bpk}
If $k\in\mathcal{K}$ and $p\in[1,+\infty]$ we denote by $B_{p,k}(\R^n)$ the set of all distributions $w\in\S'(\R^n)$ such that $\widehat{w}$ is a function and 
\[
\Vert w\Vert_{p,k}=(2\pi)^{-n}\Vert k\widehat{w}\Vert_p <	\infty.
\]
\end{definition}
The inequality \eqref{est_Hoer} says that $\wt{P_\eps}$ is a tempered weight function for each $\eps$ so it is meaningful to consider the sets $B_{\infty,\wt{P_\eps}}(\R^n)$ of distributions as we will see in the next theorem.

\begin{theorem}
\label{theo_fund_P}
To every differential operator $P(D)$ with coefficients in $\wt{\C}$ such that $\wt{P}(\xi)$ is invertible in some $\xi_0\in\R^n$ there exists a fundamental solution $E\in\Lb(\Gc(\R^n),\wt{\C})$. More precisely, to every $c>0$ and $(P_\eps)_\eps$ representative of $P$ there exists a fundamental solution $E$ given by a net of distributions $(E_\eps)_\eps$ such that $E_\eps/\cosh(c|x|)\in B_{\infty,\wt{P_\eps}}(\R^n)$ and for all $\eps$
\[
\biggl\Vert \frac{E_\eps}{\cosh(c|x|)}\biggr\Vert_{\infty,\wt{P_\eps}}\le C_0,
\]
where the constant $C_0$ depends only on $n,m$ and $c$.
\end{theorem} 

The proof of Theorem \ref{theo_fund_P} requires some technical preparation which consists in recalling and applying some classical results due to H\"ormander \cite[Chapter III]{Hoermander:63}, \cite[Chapter X]{Hoermander:V2} to the representing nets $(P_\eps)_\eps$ of polynomials with generalized constant coefficients.

\begin{lemma}
\label{lemma_Hoer_1}
Let $A$ be a bounded subset of $\R^n$ such that no polynomial of degree $\le m$ vanishes in $A$ without vanishing identically. Let us set 
\[
A'=\{k\theta/m;\, 0\le k\le m,\ \theta\in A\},
\]
where $k$ is an integer. Then there is a constant $C>0$ such that
\[
\wt{p}(\xi)\le C\sup_{\theta\in A'}\inf_{|z|=1}|p(\xi+ z\theta)|
\]
for all polynomials $p$ of degree $\le m$ and every complex $\xi$.
\end{lemma}
Lemma \ref{lemma_Hoer_1} is proved in \cite[Chapter III]{Hoermander:63} and applies to nets $(P_\eps)_\eps$ of polynomials of degree $\le m$. In particular, there exists $C>0$ such that for all nets $(P_\eps)_\eps$ the equality
\beq
\label{eq_lemma_1}
\wt{P_\eps}(\xi)\le C\sup_{\theta\in A'}\inf_{|z|=1}|P_\eps(\xi+ z\theta)|
\eeq
holds for all $\xi\in\C^n$ and for all $\eps\in(0,1]$. We are now in the position of stating the following proposition whose proof is a straightforward application of Theorems 3.1.1 and 3.1.2 in \cite{Hoermander:63} to the net $P_\eps(D)$.
\begin{proposition}
\label{lemma_Hoer_2}
Let $A'$ be a finite subset of the sphere $|\xi|<c$ such that \eqref{eq_lemma_1} is valid for every net $(P_\eps)_\eps$ of polynomials of degree $\le m$. Let us fix a net $(P_\eps)_\eps$ and let $\varphi_{\theta,\eps}$, $\theta\in A'$, $\eps\in(0,1]$ be measurable functions in $\R^n$ such that $\varphi_{\theta,\eps}\ge 0$, $\sum_{\theta\in A'}\varphi_{\theta,\eps} =1$ and 
\beq
\label{impl}
\varphi_{\theta,\eps}(\xi)>0\quad \Rightarrow \quad \wt{P_\eps}(\xi)\le C\inf_{|z|=1}|P_\eps(\xi+ z\theta)|.
\eeq
Then the formula 
\beq
\label{form_Hoer}
\check{E_\eps}(u)=(2\pi)^{-n}\sum_{\theta\in A'}\int_{\R^n}\varphi_{\theta,\eps}(\xi)\, d\xi\, \frac{1}{2\pi i}\int_{|z|=1}\frac{\widehat{u}(\xi+ z\theta)}{P_\eps(\xi+ z\theta)}\frac{dz}{z}\qquad u\in\Cinfc(\R^n),
\eeq
defines a fundamental solution $E_\eps(u):=\check{E_\eps}\ast u(0)$ of $P_\eps(D)$ such that $E_\eps /\cosh(c|x|)\in B_{\infty,\wt{P_\eps}}(\R^n)$. Moreover, there exists a constant $C_0$ depending only on $n$, $m$ and $c$ such that 
\[
\biggl\Vert \frac{E_\eps}{\cosh(c|x|)}\biggr\Vert_{\infty,\wt{P_\eps}}\le C_0
\]
for all $\eps\in(0,1]$.
\end{proposition}
The proof of Theorem \ref{theo_fund_P} is at this point just a matter of combining the invertibility of the generalized weight function $\wt{P}(\xi)$ in $\xi_0$ with the nets of fundamental solutions $(E_\eps)_\eps$ of $(P_\eps(D))_\eps$ provided by Proposition \ref{lemma_Hoer_2}.

\begin{proof}[Proof of Theorem \ref{theo_fund_P}] 
Let us fix a representative $(P_\eps)_\eps$ of $P$ and consider the corresponding net $(\check{E_\eps})_\eps\in\D'(\R^n)^{(0,1]}$ given by \eqref{form_Hoer}. It determines a basic functional $\check{E}$ in $\LL(\Gc(\R^n),\wt{\C})$. Indeed, from \eqref{impl} we have that 
\[
|\check{E_\eps}(u)|\le C(2\pi)^{-n-1}\sum_{\theta\in A'}\int_{|z|=1} \int_{\R^n}\frac{|\widehat{u}(\xi+z\theta)|}{\wt{P_\eps}(\xi)}\, d\xi\, dz
\]
for all $u\in\Cinfc(\R^n)$ and all $\eps\in(0,1]$. The invertibility of $\wt{P}$ in some point $\xi_0$ yields that the estimate
\begin{multline*}
|\check{E_\eps}(u)|\le C(2\pi)^{-n-1}\sum_{\theta\in A'}\int_{|z|=1} \int_{\R^n}\frac{|\widehat{u}(\xi+z\theta)|}{\eps^N(1+C_1|\xi_0-\xi|)^{-m}}\, d\xi\, dz\\
\le C'\eps^{-N}\sum_{\theta\in A'}\int_{|z|=1}\int_{\R^n}(1+|\xi|)^m|\widehat{u}(\xi+z\theta)|\, d\xi\, dz
\le C''\eps^{-N}\, \sup_{y\in K, |\beta|\le m+n+1}|\partial^\beta u(y)|
\end{multline*} 
is valid for all $u\in\Cinf_K(\R^n)$ when $\eps$ is small enough. 
It follows that the net $(\check{E_\eps})_\eps$ gives a basic functional $\check{E}\in \LL(\Gc(\R^n),\wt{\C})$. By construction we have that $E(u)=\check{E}\ast u(0)$ is a basic functional in $\LL(\Gc(\R^n),\wt{\C})$ such that $P(D)E=\iota_d(\delta)$. Finally, let $(E_\eps)_\eps$ be the net of distributions given by $E_\eps(u)= \check{E_\eps}\ast u(0)$ where $(\check{E_\eps})_\eps$ is defined by $(P_\eps)_\eps$ as in \eqref{form_Hoer}. This net generates the basic functional $E$ and by Proposition \ref{lemma_Hoer_2} we know that $E_\eps/\cosh(c|x|)\in B_{\infty,\wt{P_\eps}}(\R^n)$ with $\Vert {E_\eps}/{\cosh(c|x|)}\Vert_{\infty,\wt{P_\eps}}\le C_0$ for all $\eps$.  
\end{proof}
Theorem \ref{theo_fund_P} entails the following solvability result.
\begin{theorem}
\label{theom_solv_dual_1}
Let $P(D)$ be a partial differential operator with coefficients in $\wt{\C}$ such that $\wt{P}$ is invertible in some $\xi_0\in\R^n$. Then the equation
\beq
\label{eq_1_P}
P(D)u=v
\eeq
\begin{itemize}
\item[(i)] has a solution $u\in\G(\R^n)$ if $v\in\Gc(\R^n)$,
\item[(ii)] has a solution $u\in\Ginf(\R^n)$ if $v\in\Gcinf(\R^n)$,
\item[(iii)] has a solution $u\in\LL(\Gc(\R^n),\wt{\C})$ if $v\in\LL(\G(\R^n),\wt{\C})$,
\item[(iv)] has a solution $u\in\Lb(\Gc(\R^n),\wt{\C})$ if $v\in\Lb(\G(\R^n),\wt{\C})$.
\end{itemize}
\end{theorem}
\begin{proof}
Let $E\in\Lb(\Gc(\R^n),\wt{\C})$ be a fundamental solution of $P(D)$ whose existence is guaranteed by Theorem
\ref{theo_fund_P}. Combining Proposition \ref{prop_conv_1} with Proposition \ref{conv_1} we have that $u=v\ast E$ is a solution of equation \eqref{eq_1_P}. Indeed, $P(D)u=P(D)(v\ast E)=v\ast\iota_d(\delta)=v$. More precisely, $u\in\G(\R^n)$ if $v\in\Gc(\R^n)$ and $u\in\Ginf(\R^n)$ if $v\in\Gcinf(\R^n)$. When $v$ is a functional in $\LL(\G(\R^n),\wt{\C})$ then $u\in\LL(\Gc(\R^n),\wt{\C})$ and in addition $u$ is a basic functional if $v$ is basic itself.  
\end{proof}
\begin{remark}
\label{rem_HO}
The first assertion of Theorem \ref{theom_solv_dual_1} was already proven by H\"ormann and Oberguggenberger in \cite[Theorem 7.7]{HO:03}. In the course of the proof the authors define a representative $(u_\eps)_\eps$ of the solution $u$ as the convolution $(E_\eps\ast v_\eps)_\eps$, where $(v_\eps)_\eps$ is a representative of $v$ and $(E_\eps)_\eps$ a net of distributional fundamental solutions of $P_\eps(D)$. Since the dual $\LL(\Gc(\R^n),\wt{\C})$ does not appear in their mathematical framework, they do not view $(E_\eps)_\eps$ as a net defining a generalized object. They immediately consider the class in $\G(\R^n)$ generated by $(E_\eps\ast v_\eps)_\eps$.
\end{remark}
A deeper investigation of the solvability of the equation $P(D)u=v$, when the right-hand side has compact support, is postponed to the appendix at the end of the paper and modelled on the classical sources \cite{Hoermander:63, Hoermander:V2}.

\subsection{Application to evolution operators}
\label{subsec_evol}
In the sequel we set $H_n=\{x\in\R^n:\, x_n\ge 0\}$.  
\begin{definition}
\label{def_evol_op}
A partial differential operator with constant Colombeau coefficients, defined on $\R^n$, is called an evolution operator with respect to $H_n$ if it has a fundamental solution in $\Lb(\Gc(\R^n),\wt{\C})$ whose support is contained in $H_n$.
\end{definition}
The definition of evolution operator entails the following solvability results. 
\begin{theorem}
\label{theo_solv_evol}
Let $P(D)$ be an evolution operator with respect to $H_n$ .
\begin{itemize}
\item[(i)] If $v\in\Gc(\R^n)$ then the equation $P(D)u=v$ has a solution $u\in\G(\R^n)$ with $$\supp\, u\subseteq\{x\in\R^n:\, x_n\ge\inf\{y_n:\, y\in\supp\, v\}\}.$$
\item[(ii)] $(i)$ holds with $\Gc(\R^n)$ and $\G(\R^n)$ substituted by $\LL(\G(\R^n),\wt{\C})$ and $\LL(\Gc(\R^n),\wt{\C})$ respectively.
\item[(iii)] If $v\in\LL(\Gc(\R^n),\wt{\C})$ has support contained in a closed cone $\Gamma\subseteq H_n$ such that $\Gamma\cap\{x:\, x_n=0\}=\{0\}$ then the equation $P(D)u=v$ has a solution in $\LL(\Gc(\R^n),\wt{\C})$ with $\supp\, u\subseteq H_n$.
\end{itemize}
\end{theorem}
\begin{proof}
Let $E\in\Lb(\Gc(\R^n),\wt{\C})$ be a fundamental solution of $P(D)$ with $\supp\, E\subseteq H_n$. From Theorem \ref{theom_solv_dual_1} we have that $u=v\ast E$ is a solution of the equation $P(D)u=v$ which belongs to $\G(\R^n)$ when $v\in\Gc(\R^n)$ and to $\LL(\Gc(\R^n),\wt{\C})$ when $v\in\LL(\G(\R^n),\wt{\C})$. In particular by the theorem of supports we obtain that $\supp\, u\subseteq\supp\, E+\supp\, v\subseteq H_n+\supp\,v$. This means that $\supp\, u\subseteq\{x\in\R^n:\, x_n\ge\inf\{y_n:\, y\in\supp\, v\}\}$. Finally, assume that $v\in\LL(\Gc(\R^n),\wt{\C})$ has support contained in a closed cone $\Gamma\subseteq H_n$ such that $\Gamma\cap\{x:\, x_n=0\}=\{0\}$. Since $\supp\, E\subseteq H_n$ then by Proposition \ref{prop_conv_proper} and Remark \ref{rem_cone_case} we conclude that $v\ast E$ is a well-defined element of $\LL(\Gc(\R^n),\wt{\C})$ with support contained in $H_n$ and clearly a solution of $P(D)u=v$.
\end{proof}
We provide now a condition on the generalized polynomial of $P(D)$ which is sufficient to claim that $P(D)$ is an evolution operator with respect to $H_n$. In the course of the proof of Theorem \ref{theo_evol_op} we will use the fact that if $T\in\LL(\GS(\R^n),\wt{\C})$ and $u\in\G(\R^n)$ then $uT(v):=T(uv)$ defines a functional in $\LL(\Gc(\R^n),\wt{\C})$.
\begin{theorem}
\label{theo_evol_op}
Let $P(\zeta)=\sum_{|\alpha|\le m}a_\alpha\zeta^\alpha$ be a polynomial on $\C^n$ with coefficients in $\wt{\C}$ fulfilling the following conditions:
\begin{itemize}
\item[(i)] there exists a choice of representatives $(a_{\alpha,\eps})_\eps$ of $a_\alpha$,

there exist a net  $(c_\eps)_\eps\in\R^{(0,1]}$ of the form $c_\eps=c\,\omega_\eps$, $c<0$, with the property
\[
\exists c_0,a>0\, \forall\eps\in(0,1]\qquad c_0\eps^a\le\omega_\eps\le \log(1/\eps)+1
\]
and there exists a constant $\eps_0>0$ such that

\[
P_\eps(\xi',\zeta_n):=\sum_{|\alpha|\le m}a_{\alpha,\eps}(\xi',\zeta_n)^\alpha\neq 0
\]
for all $\eps\in(0,\eps_0]$, for all $\xi'\in\R^{n-1}$, for all $\zeta_{n}=\xi_n+i\eta_{n}$ with $\eta_{n}<c_\eps$;   
\item[(ii)] the coefficient of the highest power of $\zeta_n$ in $P$ is an invertible element of $\wt{\C}$.
\end{itemize}
Then $P(D)$ is an evolution operator with respect to $H_n$.
\end{theorem}
\begin{proof}[Proof of Theorem \ref{theo_evol_op}] 
By assumption $(ii)$ we can write the polynomial $P(\xi',\zeta_n)$ as 
\[
A(\zeta_n^k+\zeta_n^{k-1}P_1(\xi')+...+P_k(\xi'))
\]
where $A\in\wt{\C}$ is invertible, every $P_j$ is a polynomial in $\xi'$ with complex generalized coefficients and the integer $k$ does not exceed $m$. The representing net $(P_\eps(\xi',\zeta_n))_\eps$ satisfying condition $(i)$ is of the form 
\[
A_\eps(\zeta_n-\lambda_{1,\eps})...(\zeta_n-\lambda_{k,\eps}),
\]
where $(A_\eps)_\eps$ is a representative of $A$ and the $(\lambda_{j,\eps})_\eps$ are nets of functions of $\xi'$. In addition from $(i)$ we have that ${\rm{Im}}\,\lambda_{j,\eps}(\xi')\ge c_\eps$ for all $j=1,...,k$, for all $\xi'\in\R^{n-1}$ and for all $\eps\in(0,\eps_0]$. Combining the invertibility of $A$ with the properties of $(P_\eps(\xi',\zeta_n))_\eps$ we obtain that there exist $r_1\in\R$ and $\eps_1\in(0,1]$ such that 
\beq
\label{est_P}
|P_\eps(\xi',\zeta_n)|\ge\eps^{r_1}(c_\eps-{\rm{Im}}\,\zeta_n)^k
\eeq
for all $\xi'\in\R^{n-1}$, for all $\zeta_n$ with ${\rm{Im}}\, \zeta_n < c_\eps$ and for all $\eps\in(0,\eps_1]$. 
It follows that for fixed $\eta_{n,\eps}=c'\omega_\eps$ with $c'<c$ we have
\[
|P_\eps(\xi',\xi_n+i\eta_{n,\eps})|\ge \eps^{r_1}(c_\eps-\eta_{n,\eps})^k=\eps^{r_1}(c\omega_\eps-c'\omega_\eps)^k\ge\eps^{r_1}(c-c')^k c^k_0\eps^{ak}
\]
for all $\xi'\in\R^n$, $\xi_n\in\R$ and $\eps\in(0,\eps_1]$. Hence, the net 
\[
S_\eps(\xi',\xi_n):=\begin{cases} (P_\eps(\xi',\xi_n+i\eta_{n,\eps}))^{-1} & \eps\in(0,\eps_1]\\
0 & \eps\in(\eps_1,1]
\end{cases}
\]
defines a basic functional $S$ in $\LL(\GS(\R^n),\wt{\C})$. 

Let $\eta_n$ be the real generalized number defined by $\eta_{n,\eps}=c'\omega_\eps$. As mentioned in Section \ref{section_preliminaries} the logarithmic growth condition on $\eta_n$ yields the well-definedness of $\esp_{\eta_n}:=[(\esp^{-x_n\eta_{n,\eps}})_\eps]$ as a generalized function in $\Ginf(\R^n)$. We can now compute the inverse Fourier transform of $S$ and define the basic functional 
\[
E= \esp_{\eta_n}\cdot\mF^{-1}S
\]
of $\LL(\Gc(\R^n),\wt{\C})$. Let $u\in\Gc(\R^n)$ and $\wt{u}(x):=u(-x)$. By Proposition \ref{prop_FL} we know that $\mF\LL(u)(\cdot,\cdot+i\eta_n)$ is a generalized function in $\GS(\R^n)$. Moreover, $\mF^{-1}(\esp_{\eta_n}\wt{u})=(2\pi)^{-n}\mF\LL(u)(\cdot,\cdot+i\eta_n)$ and the definitions of $S$ and $E$ entails the equalities 
\begin{multline}
\label{E_1}
E(\wt{u})=\mF^{-1}S(\esp_{\eta_n}\wt{u})=S(\mF^{-1}(\esp_{\eta_n}\wt{u}))=S((2\pi)^{-n}\mF\LL(u)(\cdot,\cdot+i\eta_n))\\=\biggl[\biggl((2\pi)^{-n}\int_{\R^n}\frac{\mF\LL(u_\eps)(\xi',\xi_n+i\eta_{n,\eps})}{P_\eps(\xi',\xi_n+i\eta_{n,\eps})}\, d\xi'\, d\xi_n\biggr)_\eps\biggr].
\end{multline}
Note that \eqref{est_P} combined with Proposition \ref{prop_FL}$(i)$ allows to write the integral in \eqref{E_1} for all $\eps\in(0,\eps_1]$. By Cauchy's theorem this integral does not depend on the constant $c'<c$ which appears in the definition of $(\eta_{n,\eps})_\eps$. 

We now prove that $E$ is a fundamental solution of $P(D)$. From \eqref{E_1} and the equality \eqref{eq_FL} we have that 
\begin{multline}
\label{E_2}
P(D)(E)(\wt{u})=E({\,}^tP(D)\wt{u})=E((P(D)u)\wt{\,}\,)=\biggl[\biggl((2\pi)^{-n}\int_{\R^n}\frac{\mF\LL(P_\eps u_\eps)(\xi',\xi_n+i\eta_{n,\eps})}{P_\eps(\xi',\xi_n+i\eta_{n,\eps})}\, d\xi'\, d\xi_n\biggr)_\eps\biggr]\\
=\biggl[\biggl((2\pi)^{-n}\int_{\R^n}{\mF\LL(u_\eps)(\xi',\xi_n+i\eta_{n,\eps})}\, d\xi'\, d\xi_n\biggr)_\eps\biggr]=(2\pi)^{-n}\int_{\R^n}\mF\LL(u)(\xi',\xi_n+i\eta_{n})\, d\xi'\, d\xi_n.
\end{multline}
Applying again Cauchy's theorem at the level of representatives we conclude that 
\beq
\label{E_3}
(2\pi)^{-n}\int_{\R^n}\mF\LL(u)(\xi',\xi_n+i\eta_n)\, d\xi'\, d\xi_n =(2\pi)^{-n}\int_{\R^n}\mF\LL(u)(\xi',\xi_n)\, d\xi'\, d\xi_n =\int_{\R^n}\mF(u)(\xi)\, \dslash\xi =u(0).
\eeq
Finally, a combination of \eqref{E_2} with \eqref{E_3} leads to 
\[
P(D)(E)(\wt{u})=u(0)=\wt{u}(0)=\iota_d(\delta)(\wt{u}).
\]
It remains to show that the support of $E$ is contained in the region $H_n$. Let $u\in\Gc(\R^n)$ with $\supp\, u\subseteq\{x\in\R^n:\, x_n>0\}$. By Remark \ref{rem_FL} we know that there exist $\eps_2\in(0,1]$ and $N\in\N$ such that 
\[
\forall(\xi',\xi_n)\in\R^n\, \forall\eps\in(0,\eps_2]\, \forall c'<c\, \forall\eta_\eps=c'\omega_\eps\qquad\qquad |\mF\LL(u_\eps)(\xi',\xi_n+i\eta_{n,\eps})|\le \eps^{-N}(1+|\xi'|+|\xi_n|)^{-n-1}. 
\]
Hence, from \eqref{E_1} and \eqref{est_P} we conclude that $E(\wt{u})$ has a representative $(E(\wt{u}))_\eps$ satisfying the estimate
\beq
\label{final_est}
|(E(\wt{u}))_\eps|\le \eps^{-N-r_1}\frac{1}{(2\pi)^n(c-c')^k c^k_0\eps^{ak}}\int_{\R^n}\frac{1}{(1+|\xi'|+|\xi_n|)^{n+1}}\, d\xi'\, d\xi_n\le C\eps^{-N-r_1-ak}\frac{1}{(c-c')^k}
\eeq
for all $\eps\in(0,\eps_2]$ and $c'<c$. Since the left hand-side of \eqref{final_est} does not depend on $c'<c$, letting $c'$ tends to $-\infty$ we conclude that the net $(E(\wt{u}))_\eps)_\eps$ is identically $0$ when $\eps$ belongs to the interval $(0,\eps_2]$. This proves that $E(u)=0$ for all $u\in\Gc(\R^n)$ with $\supp\, u\subseteq\{x\in\R^n:\, x_n<0\}$ or in other words that $\supp\, E\subseteq H_n$.
\end{proof}

\begin{example}
We give a few examples of nets $(P_\eps)_\eps$ and corresponding generalized operators $P(D)$ which satisfy the assumptions of the previous theorem.

Let $(\omega_\eps)_\eps\in\C^{(0,1]}$ and 
\beq
\label{ex_1_op}
P_\eps(\xi',\xi_n+i\eta_n)=i\xi_n-\eta_n+\omega_\eps|\xi'|^2.
\eeq
Condition $(ii)$ is trivially fulfilled. If ${\rm{Re}}\,\omega_\eps\ge 0$ for all $\eps$ and $\eta_n<c<0$ then $P_\eps(\xi',\xi_n+i\eta_n)$ satisfies condition $(i)$. Indeed,
\[
{\rm{Re}}(P_\eps(\xi',\xi_n+i\eta_n))=\Re(\omega_\eps)|\xi'|^2-\eta_n>0
\]
for all $\eps\in(0,1]$. Note that \eqref{ex_1_op} gives the heat operator for $\omega_\eps=1$ and the Schr\"odinger operator for $\omega_\eps=i$.

The previous example can be easily generalized without loosing the evolution operator's property. It suffices to take an invertible real generalized number $a\in\wt{\R}$ in the coefficient of $\xi_n+i\eta_n$. More precisely, let 
\beq
\label{ex_2_op}
P_\eps(\xi',\xi_n+i\eta_n)=ia_\eps\xi_n-a_\eps\eta_n+\omega_\eps|\xi'|^2,
\eeq
where $(a_\eps)_\eps\in\R^{(0,1]}$ has the property $a_\eps\ge \eps^{s}$ for some $s\in\R$ and for all $\eps\in(0,1]$. Again, if ${\rm{Re}}\,\omega_\eps\ge 0$ for all $\eps$ and $\eta_n<c<0$ we obtain that ${\rm{Re}}(P_\eps(\xi',\xi_n+i\eta_n))>0$.
Note that due to the simple structure of the examples \eqref{ex_1_op} and \eqref{ex_2_op} one finds a net $(c_\eps)_\eps$ which is constant. This does not happen in the following case.

Let $P(D)=\partial_t+a\partial_x+b$, with $a,b\in\wt{\R}$. Fixing a choice of representatives $(a_\eps)_\eps$ and $(b_\eps)_\eps$ we can write
\[
P_\eps(\xi_1,\xi_2+i\eta_2)=i(\xi_1+a_\eps\xi_2)-a_\eps\eta_2+b_\eps.
\]
It is clear that in order to fulfill the second condition of Theorem \ref{theo_evol_op} we have to assume that $a$ is invertible in $\wt{\R}$. Moreover, we easily see that $P_\eps(\xi_1,\xi_2+i\eta_2)\neq 0$ for all $(\xi_1,\xi_2)$ in $\R^2$ and for all $\eps$ if and only if $-a_\eps\eta_2+b_\eps\neq 0$. Hence, we can assume $-a_\eps\eta_2+b_\eps>0$. This defines the net $\omega_\eps=\frac{b_\eps}{a_\eps}$ which satisfies the hypotheses of the theorem if for instance there exist $c,c_0,a_0>0$ such that $c_0\eps^a_0\le \frac{b_\eps}{ca_\eps}\le \log(1/\eps)+1$. This is the case of $\frac{1}{c}\le a_\eps\le c_1$ and $c_0\eps^{a_0}\le\frac{b_\eps}{cc_1}\le \frac{\log(1/\eps)+1}{cc_1}$. 
\end{example}

\section{$\G$- and $\Ginf$-hypoellipticity and ellipticity}
\label{section_hyp}
\begin{definition}
\label{def_hyp_op}
Let $P(D)$ be a partial differential operator with coefficients in $\wt{\C}$. $P(D)$ is said to be \emph{$\G$-hypoelliptic} if for any open subset $\Om$ of $\R^n$ the equality 
\beq
\label{G_hyp}
\singsupp_\G\, P(D)T=\singsupp_\G\, T
\eeq
holds for all basic functionals $T\in\LL(\Gc(\Om),\wt{\C})$. Analogously, $P(D)$ is said to be \emph{$\Ginf$-hypoelliptic} if for any open subset $\Om$ of $\R^n$, 
\beq
\label{Ginf_hyp}
\singsupp_{\Ginf}\, P(D)T=\singsupp_{\Ginf}\, T
\eeq
for all basic functionals $T\in\LL(\Gc(\Om),\wt{\C})$
\end{definition}
Definition \ref{def_hyp_op} can be equivalently stated by requiring that \eqref{G_hyp} and \eqref{Ginf_hyp} are valid for all basic functionals of $\LL(\Gc(\R^n),\wt{\C})$. From the pseudolocality property of $P(D)$ it follows that $P(D)$ is $\G$-hypoelliptic if and only if for all open subsets $X$ of $\R^n$ and all basic functionals $T$ of $\LL(\Gc(\R^n),\wt{\C})$, $P(D)T\vert_X\in\G(X)$ implies $T|_X\in\G(X)$. Analogously, replacing $\G$ with $\Ginf$ we obtain another equivalent formulation of the $\Ginf$-hypoellipticity of $P(D)$.

\subsection{Fundamental solutions of $\G$- and $\Ginf$-hypoelliptic operators}

For operators with constant Colombeau coefficients the $\G$-hypoellipticity as well as the $\Ginf$-hypoellipticity may be characterized making use of the fundamental solutions. We say that $F\in\Lb(\Gc(\R^n),\wt{\C})$ is a $\G$-parametrix of $P(D)$ if 
\[
P(D)F-\iota_d(\delta)\in\G(\R^n)
\]
and that $F\in\Lb(\Gc(\R^n),\wt{\C})$ is a $\Ginf$-parametrix of $P(D)$ if 
\[
P(D)F-\iota_d(\delta)\in\Ginf(\R^n).
\]
\begin{theorem}
\label{theo_hypo}
Let $P(D)$ be a partial differential operator with constant Colombeau coefficients such that the function $\wt{P}$ is invertible in some point of $\R^n$. The following assertions are equivalent:
\begin{itemize}
\item[(i)] the operator $P(D)$ is $\G$-hypoelliptic in $\R^n$,
\item[(ii)] the operator $P(D)$ admits a fundamental solution $E\in\Lb(\Gc(\R^n),\wt{\C})$ with $\singsupp_\G\, E\subseteq\{0\}$,
\item[(iii)] the operator $P(D)$ admits a $\G$-parametrix $F\in\Lb(\Gc(\R^n),\wt{\C})$ with $\singsupp_\G\, F\subseteq\{0\}$.
\end{itemize}
The same kind of equivalence holds with $\Ginf$-hypoelliptic, $\singsupp_{\Ginf}$ and $\Ginf$-parametrix in place of 
$\G$-hypoelliptic, $\singsupp_{\G}$ and $\G$-parametrix respectively.
\end{theorem}
\begin{proof}
We begin by considering the $\G$-case. Since $\supp_\G\, \iota_d(\delta)=\{0\}$ from Theorem \ref{theo_fund_P} we have that $(i)\Rightarrow (ii) \Rightarrow (iii)$. We now want to prove that $(iii)$ implies $(i)$. Let $F\in\Lb(\Gc(\R^n),\wt{\C})$ be a $\G$-parametrix of $P(D)$ with $\singsupp_\G\, F\subseteq\{0\}$ and $\psi\in\Cinfc(\R^n)$ be a cut-off function identically $1$ in a neighborhood of the origin. Then, $(1-\psi)F\in\G(\R^n)$ and $\psi F\in\Lb(\G(\R^n),\wt{\C})$. It follows that 
\[
P(D)F=P(D)(1-\psi)F+P(D)\psi F=\iota_d(\delta)+v,
\]
where $P(D)(1-\psi)F$ and $v$ belong to $\G(\R^n)$.\\
We have to prove that for all open subsets $X$ of $\R^n$ and all basic functionals $T$ of $\LL(\Gc(\R^n),\wt{\C})$ if $P(D)T\vert_X\in\G(X)$ then $T|_X\in\G(X)$. Equivalently we shall show that $T|_{X_1}\in\G(X_1)$ for any relatively compact open subset $X_1$ of $X$. Let $\alpha\in\Cinfc(X)$ be a cut-off function identically $1$ in a neighborhood of $\overline{X_1}$. By construction $P(D)\alpha T\in\Lb(\G(\R^n),\wt{\C})$ and $(P(D)\alpha T)|_{X_1}\in\G(X_1)$. Computing the convolution between $P(D)\alpha T$ and $F$ we have that
\[
P(D)\alpha T\ast F=\alpha T\ast P(D)F=\alpha T\ast (\iota_d(\delta)+v)=\alpha T +\alpha T\ast v
\]
and therefore
\[
\alpha T= P(D)\alpha T\ast F+w
\]
for some $w\in\G(\R^n)$. The assertion in Proposition \ref{prop_conv_1}  concerning the $\G$-singular support of the product of convolution and the properties of $F$ lead to 
\beq
\label{supp_G}
\singsupp_\G\, \alpha T\subseteq \singsupp_\G\, P(D)\alpha T + 0=\singsupp_\G\, P(D)\alpha T.
\eeq 
Since $(P(D)\alpha T)|_{X_1}\in\G(X_1)$ the inclusion \eqref{supp_G} entails $\singsupp_\G\, \alpha T\subseteq \R^n\setminus X_1$. This shows that $T|_{X_1}\in\G(X_1)$.

The chain of implications $(i)\Rightarrow (ii) \Rightarrow (iii)$ is also clear in the $\Ginf$-case. Let now $F$ be a basic functional in $\LL(\Gc(\R^n),\wt{\C})$ which is a $\Ginf$-parametrix of $P(D)$ with $\singsupp_{\Ginf}\, F\subseteq\{0\}$.
Then $P(D)F=\iota_d(\delta)+v$, where $v\in\Ginf(\R^n)$. We have to prove that for all open subsets $X$ of $\R^n$ and all basic functionals $T$ of $\LL(\Gc(\R^n),\wt{\C})$ if $P(D)T\vert_X\in\Ginf(X)$ then $T|_X\in\Ginf(X)$. Taking $X_1$ and $\alpha$ as above we conclude that $\alpha T= P(D)\alpha T\ast F+w$ for some $w\in\Ginf(\R^n)$. Hence 
\[
\singsupp_{\Ginf}\, \alpha T\subseteq \singsupp_{\Ginf}\, P(D)\alpha T\subseteq \R^n\setminus X_1
\]
which implies that $T|_{X_1}\in\Ginf(X_1)$.
\end{proof}
In the statement of the previous theorem we actually know that $\singsupp_\G\, E=\{0\}$ and $\singsupp_\G\, F=\{0\}$. Moreover, all the fundamental solutions of a $\G-$ or $\Ginf$-hypoelliptic operator have the same support properties.
\begin{corollary}
\label{coro_hypo}
Let $P(D)$ be a partial differential operator with constant Colombeau coefficients such that the function $\wt{P}$ is invertible in some point of $\R^n$.
\begin{itemize}
\item[(i)] If $P(D)$ is $\G$-hypoelliptic in $\R^n$ then all its fundamental solutions in $\Lb(\Gc(\R^n),\wt{\C})$ have the set $\{0\}$ as $\G$-singular support.
\item[(ii)] If $P(D)$ is $\Ginf$-hypoelliptic in $\R^n$ then all its fundamental solutions in $\Lb(\Gc(\R^n),\wt{\C})$ have the set $\{0\}$ as $\Ginf$-singular support.  
\end{itemize}
\end{corollary} 
\begin{proof}
By Theorem \ref{theo_hypo}$(ii)$ we know that when $P(D)$ is $\G$-hypoelliptic then it admits a fundamental solution $E\in\Lb(\Gc(\R^n),\wt{\C})$ with $\singsupp_\G\, E=\{0\}$. Let $T\in\Lb(\Gc(\R^n),\wt{\C})$ be another fundamental solution of $P(D)$. The difference $E-T$ is a solution of the homogeneous equation given by $P(D)$ then $\singsupp_\G(E-T)=\emptyset$. This means that $T=E+v$, where $v\in\G(\R^n)$ and therefore $\singsupp_\G\, T=\singsupp_\G\, E=\{0\}$. The proof of the second assertion consists in replacing $\G$- with $\Ginf$-.
\end{proof}
\subsection{$\G$- and $\Ginf$-elliptic operators}

We now deal with the important class of $\G$-elliptic and $\Ginf$-elliptic operators. Before stating the definition we recall that $r\in\wt{\C}$ is slow scale-invertible if there exists a slow scale net $(\omega_\eps)_\eps$ and a representative $(r_\eps)_\eps$ of $r$ such that $|r_\eps|\ge \omega_\eps^{-1}$ for $\eps$ small enough.
\begin{definition}
\label{def_elliptic}
A partial differential operator $P(D)$ of order $m$ with coefficients in $\wt{\C}$ is said to be $\G$-elliptic if the generalized number
\beq
\label{ellp_inf}
\biggl[\biggl(\inf_{\xi\in\R^n, |\xi|=1}|P_{m,\eps}(\xi)|\biggr)_\eps\biggr]
\eeq
is invertible.\\
It is said to be $\Ginf$-elliptic if the generalized number in \eqref{ellp_inf} is slow scale-invertible.
\end{definition}
Since, given two different representatives $(P_\eps)_\eps$ and $(P'_\eps)_\eps$ of $P$ the inequality
\[
\big|\inf_{|\xi|=1}|P_{m,\eps}(\xi)|-\inf_{|\xi|=1}|P'_{m,\eps}(\xi)|\big|\le\sup_{|\xi|=1}|(P_{m,\eps}-P'_{m,\eps})(\xi)|
\]
holds for all $\eps$, it follows that the generalized number in \eqref{ellp_inf} does not depend on the choice of the representatives of the polynom $P$ but on the operator $P(D)$. In particular, Definition \ref{def_elliptic} means that for any choice of representatives of the coefficients of $P(\xi)$ the net $(P_{m,\eps})_\eps$ satisfies the estimate
\beq
\label{est_P_G_ell}
|P_{m,\eps}(\xi)|\ge \eps^r,\qquad\qquad |\xi|=1,\, \eps\in(0,\eta]
\eeq
when $P(D)$ is $\G$-elliptic and the estimate
\beq
\label{est_P_Ginf_ell}
|P_{m,\eps}(\xi)|\ge c_\eps^{-1},\qquad\qquad |\xi|=1,\, \eps\in(0,\eta],
\eeq
with some slow scale net $(c_\eps)_\eps$, when $P(D)$ is $\Ginf$-elliptic. Note that it is not restrictive to assume that $\inf_\eps c_\eps\ge 2$.
\begin{theorem}
\label{theo_ellip}
\leavevmode
\begin{itemize}
\item[(i)] If $P(D)$ is a $\G$-elliptic operator with coefficients in $\wt{\C}$ then it is $\G$-hypoelliptic in $\R^n$.
\item[(ii)] If $P(D)$ is a $\Ginf$-elliptic operator with coefficients in $\G^\ssc_{\C}$ then it is $\Ginf$-hypoelliptic in $\R^n$.
\end{itemize}
\end{theorem}
The proof of Theorem \ref{theo_ellip} makes use of the following two lemmas.
\begin{lemma}
\label{lemma_elliptic}
\leavevmode
\begin{itemize}
\item[(i)] Let $P(D)$ be a $\G$-elliptic operator of order $m$ with coefficients in $\wt{\C}$. Then there exist $M\in\N$, $a\in\R$ and $\eta\in(0,1]$ such that
\[
|P_{\eps}(\xi)|\ge \eps^a |\xi|^m
\]
for all $\xi\in\R^n$ with $|\xi|\ge \eps^{-M}$ and for all $\eps\in(0,\eta]$.
\item[(ii)] If $P(D)$ is a $\Ginf$-elliptic operator of order $m$ with coefficients in $\G^\ssc_\C$ then there exist two slow scale nets $(\omega_\eps)_\eps$ and $(s_\eps)_\eps$ and a constant $\eta>0$ such that 
\[
|P_{\eps}(\xi)|\ge \omega_\eps^{-1} |\xi|^m
\]
for all $\xi\in\R^n$ with $|\xi|\ge s_\eps$ and for all $\eps\in(0,\eta]$.
\end{itemize}
\end{lemma}
\begin{proof}
\leavevmode
\begin{trivlist}
\item[(i)] Combining \eqref{est_P_G_ell} with the homogeneity of $P_{m,\eps}(\xi)$ we have that $|P_{m,\eps}(\xi)|\ge \eps^r|\xi|^m$ for all $\eps\in(0,\eta_1]$ with $\eta_1$ small enough and for all $\xi\in\R^n$. $P_\eps(\xi)$ can be written as $P_{m,\eps}(\xi)+P_{m-1,\eps}(\xi)$ where $P_{m-1,\eps}(\xi)=\sum_{|\alpha|\le m-1}c_{\alpha,\eps}\xi^\alpha$. The moderateness properties of the nets $(c_{\alpha,\eps})_\eps$ yield $|P_{m-1,\eps}(\xi)|\le \eps^{-N}|\xi|^{m-1}$ for all $\xi\in\R^n$ with $|\xi|\ge 1$, for some $N\in\N$ and for all $\eps\in(0,\eta_2]$. Hence for $|\xi|\ge \eps^{-N-r-1}$,
\[
|P_{m-1,\eps}(\xi)|\le \eps^{-N}\eps^{N+r+1}|\xi|^m=\eps^{r+1}|\xi|^m.
\]
It follows that for $|\xi|\ge \eps^{-N-r-1}$ and $\eps\in(0,\eta]$ with $\eta=\min(\eta_1,\eta_2,1/2)$ the estimate
\[
|P_\eps(\xi)|\ge |P_{m,\eps}|-|P_{m-1,\eps}|\ge \eps^r|\xi|^m-\eps^{r+1}|\xi|^m\ge \eps^r(1-\eps)|\xi|^m\ge \frac{\eps^r}{2}|\xi|^m\ge \eps^{r+1}|\xi|^m
\]
holds.
\item[(ii)] Analogously when $P(D)$ is $\Ginf$-elliptic from \eqref{est_P_Ginf_ell} we have that
\[
|P_{m,\eps}(\xi)|\ge c_\eps^{-1}|\xi|^m,\qquad\qquad \xi\in\R^n,\, \eps\in(0,\eta_1],
\]
and by definition of slow scale coefficient we obtain the estimate
\[
|P_{m-1,\eps}(\xi)|\le d_\eps |\xi|^{-1}|\xi|^m\le c_\eps^{-2}|\xi|^m,
\]
valid for $|\xi|\ge d_\eps c_\eps^2$ and for $\eps\in(0,\eta_2]$ with $\eta_2$ small enough. Thus, we conclude that
\[
|P_\eps(\xi)|\ge |P_{m,\eps}|-|P_{m-1,\eps}|\ge c_\eps^{-1}|\xi|^m-c_\eps^{-2}|\xi|^m=c_\eps^{-1}(1-c_\eps^{-1})|\xi|^m\ge \frac{c_\eps^{-1}}{2}|\xi|^m
\] 
for all $\xi$ with $|\xi|\ge d_\eps c_\eps^2$ and for all $\eps\in(0,\eta]$ with $\eta=\min(\eta_1,\eta_2)$.
\end{trivlist}
\end{proof}
\begin{lemma}
\label{lemma_cut}
Let $\varphi\in\Cinf(\R^n)$ such that $\varphi(\xi)=0$ for $|\xi|\le 1$ and $\varphi(\xi)=1$ for $|\xi|\ge 2$ and let $(s_\eps)_\eps$ be net of positive real numbers different from zero.
\begin{itemize}
\item[(i)] If $(s_\eps)_\eps,(s_\eps^{-1})_\eps\in\EM$ then $(\varphi(\xi/s_\eps)-1)_\eps\in\ES(\R^n)$;
\item[(ii)] If $(s_\eps)_\eps$ is a slow scale net with $\inf_\eps s_\eps>0$ then $(\varphi(\xi/s_\eps)-1)_\eps\in\ESinf(\R^n)$.
\end{itemize}
\end{lemma}
\begin{proof}
We begin by observing that $\varphi-1$ has compact support. Hence, for all $\alpha,\beta\in\N^n$ we get
\beq
\label{est_GS}
\sup_{\xi\in\R^n}|\xi^\alpha\partial^\beta(\varphi(\xi/s_\eps)-1)|\le \sup_{s_\eps\le|\xi|\le 2s_\eps}|\xi^\alpha\partial^\beta(\varphi(\xi/s_\eps)-1)|\le c\,s_\eps^{-|\beta|}(2s_\eps)^{|\alpha|}.
\eeq
The assertions $(i)$ and $(ii)$ follow easily from \eqref{est_GS}.
\end{proof}
In the proof of Theorem \ref{theo_ellip} we will refer to Lemma \ref{lemma_cut} and in particular to the fact that the net $(\varphi(\xi/s_\eps)-1)_\eps$ can define a generalized function in $\G(\R^n)$ or $\Ginf(\R^n)$ with an $\S$-moderate representative.
\begin{proof}[Proof of Theorem \ref{theo_ellip}]
$(i)$ We begin by assuming that $P(D)$ is a $\G$-elliptic operator of order $m$ with coefficients in $\wt{\C}$. By Lemma \ref{lemma_elliptic}$(i)$ we know that there exist $M\in\N$, $a\in\R$ and $\eta\in(0,1]$ such that 
\[
|P_{\eps}(\xi)|\ge \eps^a|\xi|^m
\]
for $|\xi|\ge \eps^{-M}$ and $\eps\in(0,\eta]$. It follows that taking $\varphi\in\Cinf(\R^n)$ as in Lemma \ref{lemma_cut} the net
\[
S_\eps(\xi):=\begin{cases} \frac{\varphi(\eps^M\xi)}{P_\eps(\xi)} & \eps\in(0,\eta]\\
0 & \eps\in(\eta,1]
\end{cases}
\]
determines a basic functional in $\LL(\GS(\R^n),\wt{\C})$. Therefore, ${\mF}^{-1}S\in\Lb(\GS(\R^n),\wt{\C})$ and 
$F(u):=\mF^{-1}S(u)$, $u\in\Gc(\R^n)$, is a basic functional in $\LL(\Gc(\R^n),\wt{\C})$. This is a $\G$-parametrix of $P(D)$. Indeed, the functional 
\[
P(D)F(u)-\iota_d(\delta)u=S(\mF^{-1}({\,}^tP(D)u))-\iota_d(\mF^{-1}1)(u)
\]
on $\Gc(\R^n)$ can be represented by the integral
\[
\int_{\R^n}\mF^{-1}_{\xi\to x}(\varphi(\eps^M\xi)-1)(x)u_\eps(x)\, dx.
\]
Since Lemma \ref{lemma_cut}$(i)$ implies that $v:=(\mF^{-1}_{\xi\to x}(\varphi(\eps^M\xi)-1))_\eps+\Neg(\R^n)$ is a well-defined element of $\G(\R^n)$, we conclude that $P(D)F-\iota_d(\delta)\in\G(\R^n)$. In order to complete the proof by Theorem \ref{theo_hypo}$(iii)$ it suffices to prove that $\singsupp_\G F\subseteq\{0\}$. An application of \cite[Lemma 6.8]{ChaPi:82} shows that the net of distributions $(F_\eps)_\eps:=(\mF^{-1}S_\eps)_\eps$ in $\S'(\R^n)^{(0,1]}$ which determines $F$ satisfies the following properties: for every $\alpha\in\N^n$ and every $\eps\in(0,1]$
\[
x^\alpha F_\eps\in\mathcal{C}^q(\R^n),
\]
with $q=m+|\alpha|-n-1$ and for all $K\Subset\R^n$ and every $\beta\in\N^n$ with $|\beta|\le q$ there exists $N\in\N$ such that
\[
\sup_{x\in K}|\partial^\beta(x^\alpha F_\eps)(x)|=O(\eps^{-N}).
\]
This means that away from $0$ the net $(F_\eps)_\eps$ is moderate, i.e, $(F_\eps|_{\R^n\setminus 0})_\eps\in\EM({\R^n\setminus 0})$. As a consequence $\singsupp_\G\, F\subseteq\{0\}$.

$(ii)$ When $P(D)$ is $\Ginf$-elliptic and has slow scale coefficients, we can substitute $\eps^M$ with the inverse of a slow scale net $(s_\eps)_\eps$ with $\inf_\eps s_\eps>0$ in the definition of the net $(S_\eps)_\eps$. By Lemma \ref{lemma_cut}$(ii)$ we deduce that $v:=(\mF^{-1}_{\xi\to x}(\varphi(s_\eps^{-1}\xi)-1))_\eps+\Neg(\R^n)$ belongs to $\Ginf(\R^n)$ and then $P(D)F-\iota_d(\delta)\in\Ginf(\R^n)$. Finally an inspection of the proof of Lemma 6.8 in \cite{ChaPi:82} shows that there exists a slow scale net $(\omega_\eps)_\eps$ such that
\[
\sup_{x\in K}|\partial^\beta(x^\alpha F_\eps)(x)|=O(\omega_\eps),
\]
with $\alpha$ and $\beta$ fulfilling the same assumptions as above. Hence, $(F_\eps|_{\R^n\setminus 0})_\eps\in\EMinf({\R^n\setminus 0})$ and we conclude that $\singsupp_{\Ginf}\, F\subseteq\{0\}$. 
\end{proof}
\begin{remark}
\label{rem_Mo_Gu}
For the sake of completeness we recall that the $\EMinf$-moderateness properties of the net of distributions $(F_\eps)_\eps$ have been already employed in [Proposition 5.4]\cite{HO:03} in order to investigate the $\Ginf$-regularity of the operator $P(D)$. In the proof of Theorem \ref{theo_ellip}, $(F_\eps)_\eps$ is regarded as an element of the dual $\LL(\Gc(\R^n),\wt{\C})$ and leads to $\Ginf$- as well as $\G$-regularity results for the operator $P(D)$ acting on the dual $\LL(\Gc(\R^n),\wt{\C})$. 
\end{remark}

\subsection{Necessary condition for $\G$- and $\Ginf$-hypoellipticity}

We present now a necessary condition for $\G$- and $\Ginf$-hypoellipticity. We employ an analytic method based on the estimates for the Fourier Laplace transform of a Colombeau generalized function worked out in the first section of the paper. Note that the necessary condition for $\Ginf$-hypoellipticity formulated in Theorem \ref{theom_nec_cond} is independently obtained in \cite[Theorem 5.5]{Garetto:06b} via a functional analytic method.

For technical reasons we will make use of the set of generalized points of log-type. We recall that when $\zeta\in{\wt{\C}}^n$ then ${\rm{Im}}\zeta=({\rm{Im}}\zeta_1,...,{\rm{Im}}\zeta_n)\in\wt{\R}^n$ and one can define the map $\wt{\C}^n\to\wt{\R}^n:\zeta\to{\rm{Im}}\zeta$. In particular from the section of preliminaries of this paper we have that if ${\rm{Im}}\zeta\in\wt{\R}^n$ is of log-type then $\esp^{|{\rm{Im}}\zeta|}\in\wt{\R}$ and 
$$\esp^{-ix\zeta}:=[(\esp^{-ix\zeta_\eps})_\eps]$$
is a well-defined generalized function in $\G(\R^n)$.  
\begin{theorem}
\label{theom_nec_cond}
Let $P(D)$ be a partial differential operator of order $m$ with coefficients in $\wt{\C}$ such that $\wt{P}$ is invertible in some point of $\R^n$ and let $N(P)$ the set of all zeros of $P$ in $\wt{\C}^n$ with imaginary part of log-type.
\begin{itemize}
\item[(i)] If $P(D)$ is $\G$-hypoelliptic then there exist $c\in\R$ and $a>0$ such that
\[
\val(|{\rm{Re}}\zeta|)\ge c+\val\big(\esp^{a|{\rm{Im}}\zeta|}\big)
\]
for all $\zeta\in N(P)$.
\item[(ii)] If $P(D)$ is $\Ginf$-hypoelliptic then
\[
\val(|{\rm{Re}}\zeta|)\ge 0
\]
for all $\zeta\in N(P)$.
\end{itemize}
\end{theorem}
\begin{proof}
$(i)$ By Theorem \ref{theo_hypo} we know that if $P(D)$ is $\G$-hypoelliptic then it admits a $\G$-parametrix in $\Lb(\Gc(\R^n),\wt{\C})$ which belongs to $\G$ outside the origin. Making use of a cut-off function $\psi$ identically $1$ in a neighborhood of the origin we can assume that there exists $F\in\Lb(\G(\R^n),\wt{\C})$ and $v\in\Gc(\R^n)$ such that
\[
P(D)F=\iota_d(\delta)+v
\]
in $\LL(\G(\R^n),\wt{\C})$. Let now $\zeta\in N(P)$. As observed above $\esp^{-ix\zeta}\in\G(\R^n_x)$ and therefore
\[
P(D)F(\esp^{-i\cdot\zeta})=1+v(\esp^{-i\cdot\zeta}).
\]
At the level of representatives this means that
\[
P_\eps(\zeta_\eps)\widehat{F_\eps}(\zeta_\eps)=1+\widehat{v_\eps}(\zeta_\eps)+n_\eps,
\]
where $(\zeta_\eps)$ is a representative of $\zeta$ such that $(\esp^{|{\rm{Im}}\zeta_\eps|})_\eps\in\EM$, $(n_\eps)_\eps\in\Neg$ and $\widehat{\,}$ denotes the Laplace-Fourier transform. The net of distributions $(F_\eps)_\eps\in\E'(\R^n)^{(0,1]}$ fulfills the following condition:
\[
\exists K\Subset\R^n\, \exists j\in\N\, \exists N\in\N\, \exists\eta\in(0,1]\, \forall u\in\Cinf(\R^n)\qquad |F_\eps(u)|\le \eps^{-N}\sup_{|\alpha|\le j, x\in K}|\partial^\alpha u(x)|.
\]
Hence for all $\eps\in(0,\eta]$ we obtain
\[
|\widehat{F_\eps}(\zeta_\eps)|\le \eps^{-N}(1+|\zeta_\eps|)^j\esp^{b|{\rm{Im}}(\zeta_\eps)|},
\]
where $b$ depends only on the compact set $K$. It follows that the net $(\widehat{F_\eps}(\zeta_\eps))_\eps$ is moderate and since $P(\zeta)=0$ in $\wt{\C}$ we conclude that 
\beq
\label{eq_1_v}
\widehat{v_\eps}(\zeta_\eps)=-1+n'_\eps,
\eeq
where $(n'_\eps)_\eps\in\Neg$. Assuming that the generalized function $v$ has $\supp\, v\subseteq\{x:\, |x|<a\}$ from Proposition \ref{prop_fl_gc} we have that
\[
\forall M\in\N\, \exists N\in\N\, \exists\eta\in(0,1]\, \forall\eps\in(0,\eta]\, \forall\zeta\in\C^n\qquad |\widehat{v_\eps}(\zeta)|\le \eps^{-N}(1+|\zeta|)^{-M}\esp^{a|{\rm{Im}}\zeta|}.
\]
This combined with \eqref{eq_1_v} leads to 
\[
|-1+n'_\eps|\le \eps^{-N}(1+|\zeta_\eps|)^{-1}\esp^{a|{\rm{Im}}\zeta_\eps|},
\]
where $N$ does not depend on $\zeta=[(\zeta_\eps)_\eps]\in N(P)$. Choosing $\eta$ small enough such that $|-1+n'_\eps|\ge 1/2$ for all $\eps\in(0,\eta]$ we can write
\beq
\label{final_est_G}
|\zeta_\eps|\le 2\eps^{-N}\esp^{a|{\rm{Im}}\zeta_\eps|}.
\eeq
\eqref{final_est_G} proves that there exist $c\in\R$ and $a>0$ such that 
\[
\val(|{\rm{Re}}\zeta|)\ge c+\val\big(\esp^{a|{\rm{Im}}\zeta|}\big)
\]
for all $\zeta\in N(P)$.

$(ii)$ If the operator $P(D)$ is $\Ginf$-hypoellptic then by Theorem \ref{theo_hypo} we find $F\in\Lb(\G(\R^n),\wt{\C})$ and $v\in\Gcinf(\R^n)$ such that $P(D)F=\iota_d(\delta)+v$ in $\LL(\G(\R^n),\wt{\C})$. Moreover, since $v$ is $\Ginf$-regular from Proposition \ref{prop_fl_gc} and Theorem \ref{theo_Paley_function} we obtain that 
\[
\exists N\in\N\, \forall M\in\N\, \exists\eta\in(0,1]\, \forall\eps\in(0,\eta]\, \forall\zeta\in\C^n\qquad |\widehat{v_\eps}(\zeta)|\le \eps^{-N}(1+|\zeta|)^{-M}\esp^{a|{\rm{Im}}\zeta|},
\]
with $\supp\, v\subseteq\{x:\, |x|<a\}$. Arguments analogous to the ones adopted in the first case yields that the assertion
\[
\exists N\in\N\, \forall M\in\N\, \exists\eta\in(0,1]\, \forall\eps\in(0,\eta]\qquad\qquad (1+|\zeta_\eps|)^{M}\le 2\eps^{-N}\esp^{a|{\rm{Im}}\zeta_\eps|},
\]
holds for all $\zeta\in N(P)$ with $N$ and $a$ independent of $\zeta$. Therefore for all $M\in\N$
\[
\val(|{\rm{Re}}\zeta|)\ge \frac{-N+\val\big(\esp^{a|{\rm{Im}}\zeta|}\big)}{M}
\]
or in other words $\val(|{\rm{Re}}\zeta|)\ge 0$.
\end{proof}
In the sequel we collect some remarks and examples concerning the previous statements. All the considered polynomials $P$ have the corresponding weight function $\wt{P}$ which is invertible in some point $\xi_0$.
\begin{remark}
\label{rem_nec_cond}
\leavevmode
\begin{trivlist}
\item[(i)] The assumption of log-type on ${\rm{Im}}\zeta$ cannot be dropped in the definition of the set $N(P)$. As a first example consider the operator $P(D)=-i[(\eps^r)_\eps]\partial_{x_1}+\partial_{x_2}$ where $r>0$. Its symbol is $P(\xi_1,\xi_2)=[(\eps^r)_\eps]\xi_1+i\xi_2$. By Theorem \ref{theo_ellip}$(i)$ we know that this operator is $\G$-elliptic and therefore $\G$-hypoelliptic. For any real number $c$ it is possible to find a zero $\zeta=(\zeta_1,\zeta_2)\in\wt{\C}^2$ of $P$ with ${\rm{Im}}\zeta$ not of log-type and such that $\val(|{\rm{Re}}\zeta|)<c$. Indeed, for $\zeta_1=[(\eps^{c'}+i\eps^{c'-r})_\eps]$ and $\zeta_2=[(-\eps^{c'}+i\eps^{c'+r})_\eps]$, where $c'<\min\{c,r\}$, the corresponding $\zeta\in\wt{\C}^2$ satisfies $P(\zeta)=0$, has ${\rm{Im}}\zeta=[(\eps^{c'-r})_\eps,(\eps^{c'+r})_\eps]$ which is not of log-type and 
\[
\val(|{\rm{Re}}\zeta|)=c'<c.
\]
\item[(ii)] Take now $P(D)=-i[(a_\eps)_\eps]\partial_{x_1}+\partial_{x_2}$, where $a_\eps\le a_\eps^{-1}$ and $(a_\eps^{-1})_\eps$ is a slow scale net. From Theorem \ref{theo_ellip}$(ii)$ we have that $P(D)$ is $\Ginf$-elliptic and therefore $\Ginf$-hypoelliptic. We can find a zero $\zeta\in\wt{\C}^2$ of the polynomial $P$ such that ${\rm{Im}}\zeta$ is not of log-type and $\val(|{\rm{Re}}\zeta|)<0$. In detail, $\zeta=(\zeta_1,\zeta_2)$, $\zeta_1=[(a_\eps^{-1}\eps^{-1}+i\eps^{-1})_\eps]$, $\zeta_2=[(-a_\eps\eps^{-1}+i\eps^{-1})_\eps]$ and $\val(|{\rm{Re}}\zeta|)=-1$.
\item[(iii)] The existence of a zero $\zeta$ of $P$ with classical imaginary part such that $\val(|{\rm{Re}}\zeta|)<0$ it is sufficient for deducing that the corresponding operator is not $\Ginf$-hypoelliptic. As an example take again the operator $P(D)=-i[(\eps^r)_\eps]\partial_{x_1}+\partial_{x_2}$ with $r>0$. The point $\zeta=(\zeta_1,\zeta_2)$ with $\zeta_1=[(\eps^{-r}+i)_\eps]$ and $\zeta_2=[(-\eps^{r}+i)_\eps]$ has $\val(|{\rm{Re}}\zeta|)=-r$ and the operator $P(D)$ is not $\Ginf$-hypoelliptic. Indeed, the generalized function $u=[(\esp^{(ix_1\eps^{-r}-x_2)})_\eps]$ satisfies $P(D)u=0$ but $u\not\in\Ginf(\R^2)$.
\item[(iv)] Finally, let us study the operator $P(D)=D_{x_1}-D_{x_2}$. It is not $\G$-hypoelliptic since the basic functional
\[
T(u)=\int_{\R}u(x,-x)dx
\]
in $\LL(\Gc(\R^2),\wt{\C})$ solves the equation $P(D)u=0$. We easily see that we can not control the valuation of the real part of the zeros of $P$. This is due to the fact that every $\zeta=(\zeta_1,\zeta_2)$ with ${\rm{\Re}}\zeta_1={\rm{\Re}}\zeta_2$ and ${\rm{\Im}}\zeta_1={\rm{\Im}}\zeta_2$ has the property $P(\zeta)=0$. 
\end{trivlist}
\end{remark}

\section{Examples of fundamental solutions and a structure theorem for basic functionals in the duals $\LL(\Gc(\R^n),\wt{\C})$ and $\LL(\G(\R^n),\wt{\C})$}
\label{section_ex}
This section is devoted to collect some interesting examples of fundamental solutions and to discuss their properties and applications. 
\subsection{Examples of fundamental solutions}
\subsubsection{Ordinary differential operators with Colombeau coefficients}
The simplest nontrivial ordinary differential operators (with Colombeau coefficients) are the first-order ones of the kind
\[
L=\frac{d}{dx}-a,
\]
where $a\in\wt{\C}$. We look for the fundamental solutions of $L$ that is for all functionals $T\in\LL(\Gc(\R),\wt{\C})$ such that 
\[
\frac{d}{dx}T-aT=\iota_d(\delta).
\]
We begin with the following proposition on the equation $\frac{d}{dx}T=0$ in $\LL(\Gc(\R),\wt{\C})$.
\begin{proposition}
\label{prop_d_0}
Let $T\in\LL(\Gc(\R),\wt{\C})$. If $\frac{d}{dx}T=0$ then $T=\lambda\in\wt{\C}$.
\end{proposition} 
\begin{proof}
Let $\phi_0\in\Cinfc(\R)$ such that $\int\phi_0=1$. Every $u\in\Gc(\R)$ can be written as follows:
\beq
\label{decomp_u}
u=\biggl(u-\int u(x)dx\,\phi_0\biggr)+\int u(x)dx\,\phi_0
\eeq
We denote the first and the second summand of the right-hand side of \eqref{decomp_u} by $u_1$ and $u_2$ respectively. Note that $u_1$ is the first derivative of the Colombeau generalized function
\[
v(x)=-\int^{+\infty}_x u(t)-\biggl(\int u(y)\, dy\biggr)\, \phi_0(t)\, dt,
\]
element of $\Gc(\R)$. Hence, 
\[
T(u)=T(u_1+u_2)=T(v')+T(u_2)=\int u(x)dx\, T(\phi_0).
\]
This completes the proof.
\end{proof}
\begin{proposition}
\label{prop_fund_or_1}
\leavevmode
\begin{itemize}
\item[(i)] All the fundamental solutions of the operator $L=\frac{d}{dx}$ are of the form
\[
E=\iota_d(H)+\lambda,
\]
where $\lambda\in\wt{\C}$ and $H$ is the Heaviside function.
\item[(ii)] Let $a\in\wt{\C}$ with real part of log-type. All the fundamental solutions of the operator $L_a=\frac{d}{dx}-a$ are of the form
\[
E=\iota_d(H)\esp^{ax}+\lambda\esp^{ax},
\]
with $\lambda\in\wt{\C}$.
\end{itemize}
\end{proposition}
\begin{proof}
$(i)$ It is clear that $\frac{d}{dx}\iota_d(H)=\iota_d(\delta)$. Let $F\in\LL(\Gc(\R),\wt{\C})$ be another fundamental solution of the operator $L=\frac{d}{dx}$. Then $\frac{d}{dx}(F-\iota_d(H))=0$. From Proposition \ref{prop_d_0} we have that $F-\iota_d(H)\in\wt{\C}$.

$(ii)$ Since the real part of $a$ is of log-type, $\esp^{ax}$ is a well-defined generalized function in $\G(\R)$. If $E$ is a fundamental solution of $L_a=\frac{d}{dx}-a$ then $\esp^{-ax}E$ is a fundamental solution of $L=\frac{d}{dx}$. Indeed, 
\[
\frac{d}{dx}(\esp^{-ax}E)(u)=-E(\esp^{-ax}u')=-E((\esp^{-ax}u)'+a\esp^{-ax}u)=\iota_d(\delta)(u).
\]
The first assertion of this proposition implies that $\esp^{-ax}E=\iota_d(H)+\lambda$ for some $\lambda\in\wt{\C}$ and therefore $E=\iota_d(H)\esp^{ax}+\lambda\esp^{ax}$.
\end{proof}
\begin{corollary}
\label{cor_solvab}
Let $a\in\wt{\C}$ with real part of log-type. $U=\esp^{ax}$ is the unique solution of the problem
\[
U'-aU=0,\qquad\qquad\qquad U|_{x=0}=1
\]
in $\LL(\Gc(\R),\wt{\C})$.
\end{corollary}
\begin{proof}
Clearly $\esp^{ax}$ is a solution of the problem. Assume that $T\in\LL(\Gc(\R),\wt{\C})$ is a solution as well and take a fundamental solution $E$ of the operator $L_a=\frac{d}{dx}-a$. It follows that $E+T$ is a fundamental solution of $L_a$. Proposition \ref{prop_fund_or_1}$(ii)$ yields the equality $T=\lambda \esp^{ax}$ for some $\lambda\in\wt{\C}$. Finally, since $T|_{x=0}=1$ we obtain that $\lambda=1$ and that $T=\esp^{ax}$.
\end{proof}
\begin{remark}
\label{rem_mixed_case}
The assumption of log-type behavior on $a$ gives a specific form to all the fundamental solutions of the differential operator $L_a$. This means that when ${{\rm{Re}}\,a}\in\wt{\R}$ is not of log-type we cannot exclude of finding a fundamental solution in $\Lb(\Gc (\R),\wt{\C})$ but we surely loose the freedom of generating any fundamental solutions by making $\lambda$ vary in $\wt{\C}$ as in Proposition \ref{prop_fund_or_1}$(ii)$. As an explanatory example let us consider $a\in\wt{\R}$ which is not of log-type. The net of distributions $E_\eps=\esp^{a_\eps x}H(x)-H(a_\eps)\esp^{a_\eps x}$ solves the equation $L_{a_\eps}E_\eps=\delta$ for all $\eps$ and generates a basic functional in $\LL(\Gc(\R),\wt{\C})$. Indeed for every $f\in\Cinfc(\R)$ with $\supp\, f\subseteq\{x:\, |x|\le r\}$ we have that if $a_\eps\le 0$ then
\[
|E_\eps(f)|=\biggl|\int_0^{+\infty}\esp^{a_\eps x}f(x)\, dx\biggr|\le r\sup_{|x|\le r}|f(x)|
\]
and if $a_\eps> 0$ then
\[
|E_\eps(f)|=\biggl|-\int_{-\infty}^0\esp^{a_\eps x}f(x)\, dx\biggr|\le r\sup_{|x|\le r}|f(x)|.
\]
Note that $(H(a_\eps))_\eps$ defines a generalized number in $\wt{\R}$ and makes us deal with the exponential $\esp^{a_\eps x}$ only when $a_\eps x$ is negative.
\end{remark}
Proposition \ref{prop_fund_or_1} can be extended to differential operators of higher order. As an explanatory example we consider the fundamental solutions in $\LL(\Gc(\R),\wt{\C})$ of the operator of second order
\[
L=\frac{d^2}{dx^2}+b\frac{d}{dx}+c,
\]
where $b,c\in\wt{\R}$. This requires the notions of exponential of a matrix, in particular of the matrix $M=\left(
\begin{array}{cc}
0 & 1\\
-c & -b	
\end{array}
\right)$, where the entries are generalized real numbers. We can now state the following proposition.  
We recall that $a\in\wt{\R}$ is strictly positive (i.e., positive and invertible) if there exists some representative $(a_\eps)_\eps$ and some $r>0$ such that $a_\eps\ge \eps^r$ for all $\eps$ small enough and that $a$ is strictly negative if and only if $-a$ is strictly positive. We can now state the following proposition.
\begin{proposition}
\label{prop_esp_M}
Let $b,c$ be generalized real numbers of log-type and $\Delta=b^2-4c$. The formula 
\beq
\label{matrix_exp}
\esp^{xM}:=\biggl[\biggl(\sum_{k=0}^\infty \frac{x^k}{k!}\left(
\begin{array}{cc}
0 & 1\\
-c_\eps & -b_\eps	
\end{array}
\right)^k\biggr)_\eps\biggr]
\eeq
gives a well-defined matrix of generalized functions in $\Ginf(\R)$ when  
\begin{itemize}
\item[(i)] $\Delta=0$ 
\item[or] $\Delta$ is invertible and in addition
\item[(ii)] positive 
\item[(iii)] or negative,
\end{itemize}
Moreover, 
\begin{itemize}
\item[(iii)'] when $\Delta$ is strictly negative and the log-type assumption on $c$ is dropped,
\end{itemize}
the formula \eqref{matrix_exp} defines a matrix $\esp^{xM}$ of generalized functions in $\G(\R)$.
\end{proposition}
\begin{proof} 
The proof of Proposition \ref{prop_esp_M} is essentially done by arguing at the level of representatives and applying the well-known classical results on exponentials of operators (see \cite[Chapter 3]{HS:74}). In detail, one can write $\esp^{xM}$ as
\begin{itemize}
\item[$(i)$] \[
\left(
\begin{array}{cc}
\esp^{-x\frac{b}{2}} & 0\\
0 & \esp^{-x\frac{b}{2}}	
\end{array}\right)\, \biggl(I+x\left(
\begin{array}{cc}
\frac{b}{2} & 1\\
\frac{-b^2}{4} & \frac{-b}{2}	
\end{array}\right)\biggr),
\]
\item[$(ii)$]
\[\left(
\begin{array}{cc}
1 & 1\\
\lambda_{1} & \lambda_{2}	
\end{array}\right)\, 
\left(
\begin{array}{cc}
\esp^{x\lambda_{1}} & 0\\
0& \esp^{x\lambda_{2}}	
\end{array}\right)\,
\left(
\begin{array}{cc}
\frac{\lambda_{2}}{\lambda_{2}-\lambda_{1}} & \frac{-1}{\lambda_{2}-\lambda_{1}} \\
\frac{-\lambda_{1}}{\lambda_{2}-\lambda_{1}} & \frac{1}{\lambda_{2}-\lambda_{1}}	
\end{array}\right),\qquad \lambda_{1,2}:=\frac{-b\pm\sqrt{\Delta}}{2},
\]
\item[$(iii)$]
\[
\left(
\begin{array}{cc}
0 & 1\\
\beta & \alpha	
\end{array}\right)\esp^{x\alpha}\left(
\begin{array}{cc}
\cos(\beta x) & -\sin(\beta x)\\
\sin(\beta x) & \cos(\beta x)	
\end{array}\right)\frac{-1}{\beta}\left(
\begin{array}{cc}
\alpha & -1\\
\beta & 0	
\end{array}\right),\qquad \alpha=-\frac{b}{2},\ \beta=\frac{\sqrt{-\Delta}}{2}.
\]
\end{itemize}
It is clear that if $b$ is of log-type and $c$ is an arbitrary element of $\wt{\R}$ then $\esp^{x\alpha}\in\Ginf(\R)$ and $\cos(\beta x), \sin(\beta x)\in\G(\R)$.
\end{proof}  
Under the hypotheses of Proposition \ref{prop_esp_M} the equalities $\frac{d}{dx}\esp^{xM}=M\esp^{xM}=\esp^{xM}M$ and $(\esp^{xM})^{-1}=\esp^{-xM}$ hold in the Colombeau context since they hold at the representatives' level.
\begin{remark}
\label{rem_2}
A direct application of Proposition \ref{prop_d_0} entails that if $R$ is a $n\times p$ matrix with entries in $\LL(\Gc(\R),\wt{\C})$ and $\frac{d}{dx}R=0$ then $R$ is a matrix with entries in $\wt{\C}$. In addition, a combination of Proposition \ref{prop_d_0} with Proposition \ref{prop_fund_or_1}$(i)$ proves that if $R=\left(
\begin{array}{c}
R_1\\
R_2
\end{array}\right)$, $R_1,R_2\in\LL(\Gc(\R),\wt{\C})$ is a solution of the equation $\frac{d}{dx}R=\left(
\begin{array}{c}
0\\
\iota_d(\delta)
\end{array}
\right)$
then $R_1=c_1$ and $R_2=\iota_d(H)+c_2$ with $c_1,c_2\in\wt{\C}$.
\end{remark}
It is now immediate to state Proposition \ref{prop_fund_system} whose proof is left to the reader.
\begin{proposition}
\label{prop_fund_system}
Let $M=\left(
\begin{array}{cc}
0 & 1\\
-c & -b	
\end{array}
\right)$ as in Proposition \ref{prop_esp_M}. All the solutions $U=\left(
\begin{array}{c}
u_1\\
u_2	
\end{array}
\right)$, $u_1,u_2\in\LL(\Gc(\R),\wt{\C})$ of
\beq
\label{system}
\frac{d}{dx}U=MU+\left(
\begin{array}{c}
0\\
\iota_d(\delta)
\end{array}
\right)
\eeq
are of the form
\beq
\label{form_U}
U=\iota_d(H)\esp^{xM}\left(
\begin{array}{c}
0\\
1
\end{array}\right)+\esp^{xM}\left(
\begin{array}{c}
c_1\\
c_2
\end{array}\right),
\eeq
with $c_1,c_2\in\wt{\C}$.
\end{proposition}
We finally come back to the second order operator $$L=\frac{d^2}{dx^2}+b\frac{d}{dx}+c,$$ where $b,c\in\wt{\R}$.  $E\in\LL(\Gc(\R),\wt{\C})$ is a fundamental solution of $L$ if and only if $U=\left(
\begin{array}{c}
E\\
E'
\end{array}\right)$ satisfies \eqref{system}. Under the assumption of Proposition \ref{prop_esp_M} concerning the corresponding matrix $M$, we conclude that all the fundamental solutions $E$ of $L$ can be generated as the first entry of the product of matrices in \eqref{form_U} with $c_1,c_2$ varying in $\wt{\C}$.

\subsubsection{Partial differential operators with Colombeau coefficients}
We provide some interesting example of partial differential operators with Colombeau coefficients and we investigate the $\G$- and $\Ginf$-hypoellipticity starting from a fundamental solution.
\begin{example}
\label{ex_per_CR}
Let us consider the operator
\[
P(D)=\frac{1}{a}iD_x-\frac{1}{b}D_y,
\]
where $a,b\in\wt{\R}$ are strictly positive. The basic functional $E$ given by
\[
u\to\biggl[\biggl(\frac{1}{2\pi}\int_{\R^2}\frac{a_\eps b_\eps}{a_\eps x+ib_\eps y}\,u_\eps(x,y)\, dx\, dy=\frac{1}{2\pi}\int_{0}^{2\pi}\int_0^{+\infty}(\cos\theta-i\sin\theta)\,u_\eps(a_\eps^{-1}\rho\cos\theta, b_\eps^{-1}\rho\sin\theta)\, d\rho\, d\theta\biggr)_\eps\biggr]
\]
is a fundamental solution of $P(D)$. Indeed, for $f\in\Cinf_K(\R^2)$, $K\Subset\R^2$, we can write
\begin{multline*}
\frac{1}{2\pi}\int_{\R^2}\frac{a_\eps b_\eps}{a_\eps x+ib_\eps y}\,f(x,y)\, dx\, dy\\=
\lim_{R\to 0}\frac{1}{2\pi}\int_{0}^{2\pi}\int_R^{C(K)\max\{a_\eps,b_\eps\}}(\cos\theta-i\sin\theta)f(a_\eps^{-1}\rho\cos\theta, b_\eps^{-1}\rho\sin\theta)\, d\rho\, d\theta,
\end{multline*}
where the constant $C$ depends only on the compact set $K$ and then for some $N\in\N$ and for all $\eps$ small enough we obtain the estimate
\[
\biggl|\frac{1}{2\pi}\int_{\R^2}\frac{a_\eps b_\eps}{a_\eps x+ib_\eps y}\,f(x,y)\, dx\, dy\biggr|\le {C(K)\max\{a_\eps,b_\eps\}}\sup_{(x,y)\in K}|f(x,y)|\le \eps^{-N}\sup_{(x,y)\in K}|f(x,y)|.
\]
Working at the level of representatives the action of $P(D)$ on $E$ is the following:
\begin{multline*}
\lim_{R\to 0}\frac{1}{2\pi}\int_{0}^{2\pi}\hskip-3pt\int_R^{+\infty}\hskip-10pt-\partial_\rho u_\eps(a_\eps^{-1}\rho\cos\theta, b_\eps^{-1}\rho\sin\theta)\, d\rho\, d\theta-\frac{i}{2\pi}\int_{0}^{2\pi}\hskip-3pt\int_R^{+\infty}\frac{1}{\rho}\,\partial_\theta u_\eps(a_\eps^{-1}\rho\cos\theta, b_\eps^{-1}\rho\sin\theta)\, d\rho\, d\theta\\
=\lim_{R\to 0}\frac{1}{2\pi}\int_{0}^{2\pi}u_\eps(a_\eps^{-1}R\cos\theta, b_\eps^{-1}R\sin\theta)\, d\theta= u_\eps(0,0).
\end{multline*}

Outside the origin the functional $E$ belongs to $\G$. Indeed $E|_{\R^2\setminus 0}=\frac{ab}{2\pi(ax+iby)}$. Since $\wt{P}^2(\xi_1,\xi_2)=\frac{1}{a^2}\xi_1^2+\frac{1}{b^2}{\xi_2^2}+\frac{1}{a^2}+\frac{1}{b^2}$ is invertible in every point of $\R^2$, by Theorem \ref{theo_hypo} we conclude that $P(D)$ is $\G$-hypoelliptic. Hence, every fundamental solution of $P(D)$ in $\Lb(\Gc(\R^2),\wt{\C})$ is of the form
\[
T=E+v,
\]
where $v\in\G(\R^2)$ is a solution of the homogeneous equation. More precisely, since the holomorphic generalized functions of $\G_H(\R^2)$ are defined as the solutions in $\G(\R^2)$ of the equation $$\frac{\partial u}{\partial x}+{i}\frac{\partial u}{\partial y}=0,$$ it follows that $v=w(a\cdot,b\cdot)$ with $w\in\G_H(\R^2)$.

If, in addition to the previous hypotheses, $a$ and $b$ are of slow scale type and slow scale invertible then $E|_{\R^2\setminus 0}\in\Ginf(\R^2\setminus 0)$. Hence, $P(D)$ is $\Ginf$-hypoelliptic and every fundamental solution of $P(D)$ in $\Lb(\Gc(\R^2),\wt{\C})$ is of the form $T=E+w(a\cdot,b\cdot)$, with $w\in\G_H(\R^2)$.
\end{example}
\begin{example}
\label{ex_Laplace}
We study the perturbation of the Laplace operator in $\R^2$ given by
\[
P(D)=-a_1D_{x_1}^2-a_2D_{x_2}^2,
\]
where $a_1,a_2\in\wt{\R}$ are strictly positive. The basic functional $E$ generated by the net of distributions
\[
E_\eps=\frac{1}{2\pi}\log\biggl(\sqrt{\frac{x_1^2}{a_{1,\eps}}+\frac{x_2^2}{a_{2,\eps}}}\biggr)\frac{1}{\sqrt{a_{1,\eps}}}\frac{1}{\sqrt{a_{2,\eps}}}.
\]
is a fundamental solution in $\LL(\Gc(\R^2),\wt{\C})$ of $P(D)$. We can easily check that $(E_\eps)_\eps$ defines a basic functional. Indeed, for $f\in\Cinf_K(\R^2)$, $K\Subset\R^2$, we can write the previous integral as 
\[
\lim_{R\to 0}\frac{1}{2\pi}\int_{0}^{2\pi}\int_R^{C(K)\max\{\frac{1}{\sqrt{a_{1,\eps}}},\frac{1}{\sqrt{a_{2,\eps}}}\}}\rho\log\rho\, f(\sqrt{a_{1,\eps}}\rho\cos\theta,\sqrt{a_{2,\eps}}\rho\sin\theta)\, d\rho\, d\theta,
\]
where the constant $C$ depends only on the compact set $K$. Hence, denoting $C(K)\max\{\frac{1}{\sqrt{a_{1,\eps}}},\frac{1}{\sqrt{a_{2,\eps}}}\}$ by $C_\eps(K)$ we obtain the following estimate  
\[
|E_\eps(f)|\le \big|\big(\frac{1}{2}C^2_\eps(K)\log C_\eps(K)-\frac{1}{4}{C^2_\eps(K)}\big)\big|\sup_{(x,y)\in K}|f(x,y)|\le \eps^{-N}\sup_{(x,y)\in K}|f(x,y)|,
\]
valid for some $N\in\N$ and for all $\eps$ small enough. 

The polynomial $\wt{P}^2(\xi_1,\xi_2)=(a_1\xi_1^2+a_2\xi_2^2)^2+4a_1\xi_1^2+4a_2\xi_2^2+4a_1^2+4a_2^2$ is invertible in any point $(\xi_1,\xi_2)$ of $\R^2$ and it is clear that outside the origin $E$ belongs to $\G$. Hence by Theorem \ref{theo_hypo} the operator $P(D)$ is $\G$-hypoelliptic. Moreover, if the coefficients $a_1$ and $a_2$ are of slow scale type and slow scale invertible then $E|_{\R^2\setminus 0}\in\Ginf(\R^2\setminus 0)$ and $P(D)$ is $\Ginf$-hypoelliptic. 
\end{example}

\subsection{Structure theorems for $\Lb(\Gc(\R^n),\wt{\C})$ and $\Lb(\G(\R^n),\wt{\C})$}
It is clear that a distributional fundamental solution $w$ of a classical partial differential operator $P(D)$ (regarded as a generalized operator) is a fundamental solution in the dual $\LL(\Gc(\R^n),\wt{\C})$ in the sense that $P(D)\iota_d(w)=\iota_d(\delta)$. In this subsection we investigate the structural properties of the spaces  $\Lb(\Gc(\R^n),\wt{\C})$ and $\Lb(\G(\R^n),\wt{\C})$ by making use of the distributional fundamental solution 
\[
E_k:=\frac{(x_1)_+^{k-1}...(x_n)_+^{k-1}}{(k-1)!)^n},\qquad x\in\R^n,\quad x_+:=xH(x)
\]
of the operator $(\partial_1...\partial_n)^k$. From what said above it follows that  
\[
(\partial_1...\partial_n)^k\iota_d(E_k)=\iota_d(\delta)
\]
in $\LL(\Gc(\R^n),\wt{\C})$.  

As a preliminary step to our structure investigation we introduce the notion of finite order in the dual $\LL(\Gc(\Om),\wt{\C})$. This employs the following spaces of generalized functions obtained by equipping $\G(\Om)$ and $\Gc(\Om)$ with different topologies where we fix the order of derivatives. In detail, let us fix $m\in\N$. We denote by $\G^m(\Om)$ the algebra $\G(\Om)$ equipped with the family of ultra-pseudo-seminorms $\{\mP_{K,m}\}_{K\Subset\Om}$ and by $\G_K^m(\Om)$ the space $\G_K(\Om)$ endowed with the topology of the ultra-pseudo-seminorm $\mP_{\G_K(\Om),m}$. Finally $\Gc^m(\Om)$ is the strict inductive limit of the locally convex topological $\wt{\C}$-modules $\{\G_K^m(\Om)\}_{K\Subset\Om}$. By construction it is clear that any of the previous spaces of order $m+1$ is continuously embedded in the corresponding of order $m$. It follows that $\LL(\Gc^m(\Om),\wt{\C})\subseteq\LL(\Gc^{m+1}(\Om),\wt{\C})$ and $\LL(\G^m(\Om),\wt{\C})\subseteq\LL(\G^{m+1}(\Om),\wt{\C})$. In particular, since $\Gc(\Om)\subseteq\Gc^m(\Om)$ and $\G(\Om)\subseteq\G^m(\Om)$ for all $m$ we have that $\LL(\Gc^m(\Om),\wt{\C})\subseteq\LL(\Gc(\Om),\wt{\C})$ and  
$\LL(\G^m(\Om),\wt{\C})\subseteq\LL(\G(\Om),\wt{\C})$. This means that the duals of $\Gc^m(\Om)$ and $\G^m(\Om)$ can be regarded as subspaces of $\LL(\Gc(\Om),\wt{\C})$ and $\LL(\G(\Om),\wt{\C})$ respectively.
\begin{definition}
\label{finite_order}
We call the elements of $\cup_{m\in\N}\LL(\Gc^m(\Om),\wt{\C})\subseteq\LL(\Gc(\Om),\wt{\C})$ \emph{functionals of finite order}.
\end{definition}
By the definition of $\LL(\Gc^m(\Om),\wt{\C})$ we easily see that if $T\in\LL(\Gc^m(\Om),\wt{\C})$ then $\partial^\alpha T\in\LL(\Gc^{m+|\alpha|}(\Om),\wt{\C})$.

\begin{proposition}
\label{prop_finite}
\leavevmode
\begin{itemize}
\item[(i)] Every functional in $\LL(\Gc(\Om),\wt{\C})$ with compact support is of finite order.
\item[(ii)] If $T$ is a basic functional in $\LL(\Gc(\Om),\wt{\C})$ defined by a net distributions $(T_\eps)_\eps\in\D'^m(\Om)$ such that 
\begin{multline}
\label{basic_m}
\forall K\Subset\Om\, \exists N\in\N\, \exists c>0\, \exists \eta\in(0,1]\, \forall f\in\Cinf_K(\Om)\, \forall\eps\in(0,\eta]\\ |T_\eps(f)|\le c\eps^{-N}\sup_{x\in K,\, |\alpha|\le m}|\partial^\alpha f(x)|
\end{multline}
then $T\in\LL(\Gc^m(\Om),\wt{\C})$.
\end{itemize}
\end{proposition}
\begin{proof}
$(i)$ If $T\in\LL(\Gc(\Om),\wt{\C})$ has compact support then taking a cut-off function $\psi$ identically $1$ on a neighborhood of $\supp\, T$ we can write $T(u)=\psi T(u)$ for all $u\in\Gc(\Om)$. Since $\supp\, \psi u\subseteq \supp\, \psi=K$ by the continuity of $T$ it follows that there exist $m\in\N$ and $C>0$ such that the inequality
\[
|T(u)|_\esp\le C\mP_{\G_K(\Om),m}(\psi u)
\]
holds for all $u\in\Gc(\Om)$. Hence, $T\in\LL(\Gc^m(\Om),\wt{\C})$.

$(ii)$ The assertion \eqref{basic_m} implies 
\[ 
|T(u)|_\esp\le C\mP_{\G_K(\Om),m}(u)
\]
for all $u\in\G_K(\Om)$ and for all $K\Subset\Om$. This means that $T|_{\G^m_K(\Om)}$ is continuous for all $K\Subset\Om$. Therefore, by the notion of strict inductive limit topology we conclude that $T$ is a continuous functional on $\Gc^m(\Om)$. 
\end{proof}


As a special example of functional of finite order defined as in Proposition \ref{prop_finite}$(ii)$ we have the functionals $T\in\LL(\Gc(\Om,\wt{\C})$ determinated by a moderate net $(T_\eps)_\eps$ of continuous functions, which are therefore elements of $\LL(\Gc^0(\Om),\wt{\C})$.

We are now ready to prove the following structure theorem for basic functionals.  
\begin{theorem}
\label{structure_theorem}
\leavevmode
\begin{itemize}
\item[(i)] The restriction of a basic functional in $\LL(\Gc(\R^n),\wt{\C})$ to a bounded open set $\Om\subseteq\R^n$ is a derivative of finite order of a functional in $\LL(\Gc^0(\Om),\wt{\C})$ defined by a net in $\M_{\mC(\Om)}$.
\item[(ii)] If $T\in\Lb(\G(\R^n),\wt{\C})$ then there is an integer $m\ge 0$ and a set of functionals $T_\alpha\in\LL(\Gc^0(\R^n),\wt{\C})$ defined by nets in $\M_{\mC(\R^n)}$ for $|\alpha|\le m$ such that 
\[
 T=\sum_{|\alpha|\le m}\partial^\alpha T_\alpha.
\]
\end{itemize}
\end{theorem}
\begin{proof}
$(i)$ Since $\Om$ is bounded we can find a cut-off function $\psi\in\Cinfc(\R^n)$ identically one on $\Om$. Then, $\psi T|_\Om=T|_\Om$ and $\psi T\in\Lb(\Gc(\R^n),\wt{\C})$. Let $(T_\eps)_\eps\in\M(\Cinfc(\R^n),\C)$ be a defining net of $T$. It follows that $(\psi T_\eps)_\eps$ is a net of distributions of finite order and more precisely
\beq
\label{est_psiT}
\exists N,M\in\N\, \exists C>0\, \exists\eta\in(0,1]\, \forall f\in\Cinf(\R^n)\, \forall\eps\in(0,\eta]\qquad |(\psi T_\eps)(f)|\le C\eps^{-M}\sup_{x\in\supp\psi, |\alpha|\le N}|\partial^\alpha f(x)|.
\eeq
By the considerations at the beginning of this subsection we have that
\beq
\label{net_cont}
\psi T_\eps=(\partial_1...\partial_n)^{N+2}E_{N+2}\ast \psi T_\eps
\eeq
for all $\eps\in(0,1]$. By \cite[Theorem 5.4.1]{FJ:98} we already know that the $(E_{N+2}\ast \psi T_\eps)_\eps$ in \eqref{net_cont} is a net of continuous functions. So, the theorem will follow once it is shown that $E_{N+2}\ast \psi T_\eps$ is ${\mC(\R^n)}$-moderate. This is clear by the fact that $E_{N+2}\in\mC^N(\R^n)$ and that \eqref{est_psiT} holds for all $f\in\mC^N(\R^n)$. Hence, for all $\eps$ small enough we can write
\[
\sup_{x\in K\Subset\R^n}|(E_{N+2}\ast \psi T_\eps)(x)|=\sup_{x\in K}|\psi T_\eps(E_{N+2}(x-\cdot))|\le C\eps^{-M}\sup_{x\in K,\,y\in\supp\, \psi,\, |\alpha|\le N}|\partial^\alpha_yE_{N+2}(x-y)|.
\]
Concluding, $E_{N+2}\ast\psi T|_\Om$ is a functional in $\LL(\Gc^0(\Om),\wt{\C})$ and is defined by a net in $\M_{\mC(\Om)}$.

$(ii)$ If $T$ has compact support then it can be written as $\psi T$ with $\psi\in\Cinfc(\Om)$ identically $1$ in a neighborhood of $\supp\, T$ and $\Om$ bounded subset of $\R^n$. For all $u\in\G(\R^n)$ we can write $T(u)=\psi T(u)=T(\psi u)=T|_\Om (\psi u)$. By the previous assertion we know that $T|_\Om=\partial^\alpha F$ where $F\in\LL(\Gc^0(\Om),\wt{\C})$ is defined by a net in $\M_{\mC(\Om)}$. Assume that $|\alpha|=m$. By Leibniz's theorem we get
\[
T(u)=\partial^\alpha F(\psi u)=\sum_{\alpha'\le\alpha}\partial^{\alpha'}((-1)^{|\alpha-\alpha'|}F\partial^{\alpha-\alpha'}\psi)(u),
\]  
where every $F\partial^{\alpha-\alpha'}\psi$ is a functional in $\LL(\Gc^0(\R^n),\wt{\C})$ determined by a net in $\M_{\mC(\R^n)}$.
\end{proof}
\begin{remark}
\label{rem_case_om}
It is clear that the assertions of Theorem \ref{structure_theorem} hold for basic functionals in $\LL(\Gc(\Om),\wt{\C})$. More precisely we have that the restriction of a basic functional in $\LL(\Gc(\Om),\wt{\C})$ to a relatively compact open subset $\Om'$ of $\Om$ is a derivative of finite order of a functional in $\LL(\Gc^0(\Om'),\wt{\C})$ defined by a net in $\M_{\mC(\Om')}$. From this result it follows that the statement concerning functionals with compact support is valid with $\R^n$ substituted by $\Om$.
\end{remark}

\section{Appendix:\\ solvability of the equation $P(D)u=v$ when $v$ is a basic functional in $\LL(\G(\R^n),\wt{\C})$}
\setcounter{section}{1}
\renewcommand{\thesection}{\Alph{section}}
\renewcommand{\theequation}{A.\arabic{equation}}
\setcounter{equation}{0}

\label{section_solvability}
The problem of the solvability of the equation
\beq
\label{appen_eq}
P(D)u=v
\eeq
when $v$ is a basic functional with compact support has been already approached in Section \ref{section_fund} (Theorem \ref{theom_solv_dual_1}) as a straightforward application of the existence of a fundamental solution in $\Lb(\Gc(\R^n),\wt{\C})$ for the operator $P(D)=\sum_{|\alpha|\le m}c_\alpha D^\alpha$ with constant Colombeau coefficients. In this appendix we provide a deeper investigation of the equation \eqref{appen_eq} which will involve estimate of $B_{p,k}$-type. Our results are modelled on the classical theory of $B_{p,k}$ spaces developed by H\"ormander in  \cite[Chapters II, III]{Hoermander:63}, \cite[Chapter X]{Hoermander:V2} and will give more precise analytic information on the basic functionals which solve \eqref{appen_eq} in the dual $\LL(\Gc(\R^n),\wt{\C})$.

We recall that $\mathcal{K}$ is the set of tempered weight functions defined in \cite[Definition 2.1.1]{Hoermander:63}. We begin by introducing the following notion of moderateness for nets of tempered distributions.
\begin{definition}
\label{mod_pk}
Let $k\in\mathcal{K}$, $1\le p\le\infty$ and $(P_\eps)_\eps$ be a net of polynomials of degree $m$. We say that $(T_\eps)_\eps\in\S'(\R^n)^{(0,1]}$ is $(B_{p,k\wt{P_\eps}})_\eps$-moderate if $(\widehat{T_\eps})_\eps$ is a net of functions and there exists $b\in\R$ such that $\Vert k\wt{P_\eps}\widehat{T_\eps}\Vert_p=O(\eps^b)$ as $\eps\to 0$.
\end{definition}
It is clear that when $\wt{P_\eps}$ is identically 1 then we have the notion of $B_{p,k}$-moderateness, i.e., $(T_\eps)_\eps$ is an element of $\M_{B_{p,k}(\R^n)}$. Moreover, since for any $f\in\Cinfc(\R^n)$ we can write
\[
|T_\eps(f)|\le\Vert k\wt{P_\eps}\widehat{T_\eps}\Vert_p\, \Vert{(f)^\vee}/{k\wt{P_\eps}}\Vert_{p'},\qquad\quad 1/{p}+{1}/{p'}=1
\]
from \eqref{est_inv} and the properties of the weight function $k$ it follows that any $(B_{p,k\wt{P_\eps}})_\eps$-moderate net $(T_\eps)_\eps$ defines a basic functional in $\LL(\Gc(\R^n),\wt{\C})$.

We collect now some properties of $(B_{p,k\wt{P_\eps}})_\eps$-moderate nets which will be useful in the sequel.
\begin{proposition}
\label{prop_appen}
Let $P$ be a polynomial with coefficients in $\wt{\C}$ such that $\wt{P}(\xi)$ is invertible in some point $\xi_0$ and let $(P_\eps)_\eps$ be a representative of $P$.
\begin{itemize}
\item[(i)] If $(T_\eps)_\eps\in\S'(\R^n)^{(0,1]}$ is $(B_{p,k\wt{P_\eps}})_\eps$-moderate and $\varphi\in\S(\R^n)$ then $(\varphi T_\eps)_\eps$ is $(B_{p,k\wt{P_\eps}})_\eps$-moderate. 
\item[(ii)] If $(T_{1,\eps})_\eps\in\E'(\R^n)^{(0,1]}$, with $\supp\, T_{1,\eps}\subseteq K\Subset\R^n$ for all $\eps$, is $(B_{p,k_1\wt{P_\eps}})_\eps$-moderate and $(T_{2,\eps})_\eps\in\S'(\R^n)^{(0,1]}$ is $B_{\infty,k_2}$-moderate then $(T_{1,\eps}\ast T_{2,\eps})_\eps$ is $(B_{p,k_1k_2\wt{P_\eps}})_\eps$-moderate.
\item[(iii)] Assertion $(i)$ holds when $(T_{1,\eps})_\eps$ is $(B_{p,k_1})_\eps$-moderate and $(T_{2,\eps})_\eps$ is  $(B_{\infty,k_2\wt{P_\eps}})_\eps$-moderate.
\item[(iv)] If $(T_\eps)_\eps$ is $(B_{p,k\wt{P_\eps}})_\eps$-moderate then $(P_\eps(D)T_\eps)_\eps$ is $B_{p,k}$-moderate.
\end{itemize}
\end{proposition}
\begin{proof}
$(i)$ Applying Theorem 2.2.5 in \cite{Hoermander:63} and in particular the inequality (2.2.9) to $(T_\eps)_\eps$ and $\varphi$ for fixed $\eps$ we obtain that $(u_\eps T_\eps)_\eps$ is a net of distributions in $\mB_{p,k\wt{P_\eps}}(\R^n)$ such that
\beq
\label{clas_est}
\Vert \varphi T_\eps\Vert_{p,k\wt{P_\eps}}\le (2\pi)^{-n}\Vert \varphi\Vert_{1,M_{k\wt{P_\eps}}}\, \Vert T_\eps\Vert_{p,k\wt{P_\eps}},
\eeq
where
\[
M_{k,\wt{P_\eps}}(\xi):=\sup_{\eta\in\R^n}\frac{k\wt{P_\eps}(\xi+\eta)}{k\wt{P_\eps}(\eta)}.
\]
The estimates \eqref{est_Hoer} and \eqref{est_inv} imply for $M_{k,\wt{P_\eps}}$ the bound $M_{k,\wt{P_\eps}}(\xi)\le (1+C_1|\xi|)^{m_1}$ valid for some constants $C_1$ and $m_1$, for all values of $\xi$ and for $\eps$ small enough. Making use of this result we conclude that there exist $\eps_0\in(0,1]$ such that the estimate
\[
\Vert \varphi\Vert_{1,M_{k,\wt{P_\eps}}}=\Vert M_{k\wt{P_\eps}}(\xi)\widehat{\varphi}(\xi)\Vert_1\le \int_{\R^n}(1+C_1|\xi|)^{m_1}(1+|\xi|)^{-n-1-m_1}d\xi\, \sup_{\xi\in\R^n}(1+|\xi|)^{n+1+m_1}|\widehat{\varphi}(\xi)|\le C_2
\]
holds for all $\eps\in(0,\eps_0]$. As a consequence, since by assumption $\Vert T_\eps\Vert_{p,k\wt{P_\eps}}=O(\eps^{-b})$ for some $b\in\R$ by \eqref{clas_est} we are lead to $\Vert \varphi T_\eps\Vert_{p,k\wt{P_\eps}}=O(\eps^{-b})$. This proves that $(\varphi T_\eps)_\eps$ is $(B_{p,k\wt{P_\eps}})_\eps$-moderate.

$(ii)$ We begin by observing that if the net $(T_{1,\eps})_\eps\in\E'(\R^n)^{(0,1]}$ fulfills the property $\supp\, T_{1,\eps}\subseteq K\Subset\R^n$ for all $\eps$, then it coincides with $\varphi T_{1,\eps}$ when $\varphi\in\Cinfc(\R^n)$ is a cut-off function identically $1$ in a neighborhood of $K$. Since $(T_{1,\eps})_\eps$ is  $(B_{p,k_1\wt{P_\eps}})_\eps$-moderate from the first assertion of this proposition we have that $(\varphi T_{1,\eps})_\eps$ is $(B_{p,k_1\wt{P_\eps}})_\eps$-moderate as well. Since $(\varphi T_{1,\eps})_\eps$ is a net of distributions with compact support which belongs to $B_{p,k_1\wt{P_\eps}}(\R^n)$ and $(T_{2,\eps})_\eps\in\S'(\R^n)^{(0,1]}$ is $B_{\infty,k_2}$-moderate, from Theorem 2.2.6 in \cite{Hoermander:63} and in particular the inequality (2.2.11) we have that $(\varphi T_{1,\eps}\ast T_{2,\eps})_\eps\in(B_{p,k_1k_2\wt{P_\eps}}(\R^n))^{(0,1]}$ and
\beq
\label{est_conv}
\Vert\varphi T_{1,\eps}\ast T_{2,\eps}\Vert_{p,k_1k_2\wt{P_\eps}}\le \Vert \varphi T_{1,\eps}\Vert_{p,k_1\wt{P_\eps}}\,  \Vert T_{2,\eps}\Vert_{\infty,k_2}
\eeq
for all $\eps\in(0,1]$. It follows that $(T_{1,\eps}\ast T_{2,\eps})_\eps$ is $(B_{p,k_1k_2\wt{P_\eps}})_\eps$-moderate.

$(iv)$ If $(T_\eps)_\eps$ is $(B_{p,k\wt{P_\eps}})_\eps$-moderate then $\Vert k\wt{P_\eps}\widehat{T_\eps}\Vert_p=O(\eps^b)$ for some $b\in\R$. Take now the net of tempered distributions $(P_\eps(D)T_\eps)_\eps$. Since $\widehat{P_\eps(D)T_\eps}=P_\eps\widehat{T_\eps}$ and
\[
\Vert k\widehat{P_\eps(D)T_\eps}\Vert_p=\Vert kP_\eps\widehat{T_\eps}\Vert_p\le \Vert k\wt{P_\eps}\widehat{T_\eps}\Vert_p
\]
for all $\eps\in(0,1]$, we conclude that $(P_\eps(D)T_\eps)_\eps$ is $B_{p,k}$-moderate.
\end{proof}
The notion of moderateness with respect to the net $(B_{p,k\wt{P_\eps}})$ can be expressed locally as follows.
\begin{definition}
\label{loc_mod_pk}
Let $k\in\mathcal{K}$, $1\le p\le\infty$ and $(P_\eps)_\eps$ be a net of polynomials of degree $m$. We say that $(T_\eps)_\eps\in\D'(\R^n)^{(0,1]}$ is locally $(B_{p,k\wt{P_\eps}})_\eps$-moderate (or $(T_\eps)_\eps$ is $(B^{\rm{loc}}_{p,k\wt{P_\eps}})_\eps$-moderate) if for all $\varphi\in\Cinfc(\R^n)$ the net $(\varphi{T_\eps})_\eps$ is $(B_{p,k\wt{P_\eps}})_\eps$-moderate.  
\end{definition}
Any $(B^{\rm{loc}}_{p,k\wt{P_\eps}})_\eps$-moderate net $(T_\eps)_\eps$ generates a basic functional in $\LL(\Gc(\R^n),\wt{\C})$. The results on the convolution product of Proposition \ref{prop_appen} can be extended to locally $(B_{p,k\wt{P_\eps}})_\eps$-moderate nets.
\begin{proposition}
\label{prop_appen_loc}
Under the assumption of Proposition \ref{prop_appen}, if $(T_{1,\eps})_\eps\in\E'(\R^n)^{(0,1]}$, with $\supp\, T_{1,\eps}\subseteq K\Subset\R^n$ for all $\eps$, is $(B_{p,k_1})_\eps$-moderate and $(T_{2,\eps})_\eps\in\D'(\R^n)^{(0,1]}$ is $(B^{\rm{loc}}_{\infty,k_2\wt{P_\eps}})_\eps$-moderate then $(T_{1,\eps}\ast T_{2,\eps})_\eps$ is $(B^{\rm{loc}}_{p,k_1k_2\wt{P_\eps}})_\eps$-moderate.
\end{proposition}
\begin{proof}
Let $\varphi\in\Cinfc(\R^n)$. We choose $\psi_1\in\Cinfc(\R^n)$ identically 1 in a neighborhood of the compact set $K$ and $\psi_2\in\Cinfc(\R^n)$ identically 1 in a neighborhood of $\supp\, \varphi-\supp\, \psi_1$. We can write
$\varphi(T_{1,\eps}\ast T_{2,\eps})$ as $\varphi(\psi_1 T_{1,\eps}\ast \psi_2 T_{2,\eps})$. Hence from Proposition \ref{prop_appen} we have
that $(\psi_1 T_{1,\eps})_\eps$ is $B_{p,k_1}$-moderate, $(\psi_2 T_{2,\eps})_\eps$ is $(B_{\infty,k_2\wt{P_\eps}})_\eps$-moderate and
finally $(\varphi(\psi_1 T_{1,\eps}\ast \psi_2 T_{2,\eps}))_\eps$ is $(B_{p,k_1k_2\wt{P_\eps}})_\eps$-moderate.

\end{proof}

\begin{remark}
\label{rem_fund_P}
We are now able to give a more precise description of the fundamental solution $E$ of $P(D)$ provided by Theorem \ref{theo_fund_P}. Let $P(D)$ be a partial differential operator with coefficients in $\wt{\C}$ such that $\wt{P}(\xi)$ is invertible in some $\xi_0\in\R^n$. For every representative $(P_\eps)_\eps$ of $P$ and every $c>0$ there exists a fundamental solution $E\in\Lb(\Gc(\R^n),\wt{\C})$ of $P(D)$ which is defined by a $(B^{\rm{loc}}_{\infty,\wt{P_\eps}})_\eps$-moderate net of distributions $(E_\eps)_\eps$ and such that $(E_\eps/\cosh(c|x|))_\eps$ is $(B_{\infty,\wt{P_\eps}})_\eps$-moderate.
\end{remark}
The existence of a fundamental solution with the previous moderateness properties entails the following result of solvability.
\begin{theorem}
\label{theo_solv_dual}
Let $P(D)$ be a differential operator $P(D)$ with coefficients in $\wt{\C}$ such that $\wt{P}(\xi)$ is invertible in some $\xi_0\in\R^n$ and let $v\in\Lb(\G(\R^n),\wt{\C})$. If $v$ is defined by a $B_{p,k}$-moderate net $(v_\eps)_\eps$, then the equation
\[
P(D)u=v
\]
has a solution $u\in\Lb(\Gc(\R^n),\wt{\C})$ which is given by a $(B^{{\rm{loc}}}_{p,k\wt{P_\eps}})_\eps$-moderate net $(u_\eps)_\eps$.
\end{theorem}
\begin{proof}
By Remark \ref{rem_fund_P} we know that the operator $P(D)$ admits a fundamental solution $E\in\Lb(\Gc(\R^n),\wt{\C})$ which is determined by a $({B}^{{\rm{loc}}}_{\infty,\wt{P_\eps}})_\eps$-moderate net $(E_\eps)_\eps$. Since $(v_\eps)_\eps$ is $B_{p,k}$-moderate and by Proposition \ref{prop_appen}$(i)$ it is not restrictive to assume that $\supp\, v_\eps\subseteq K\Subset\R^n$ for all $\eps$, by Proposition \ref{prop_appen_loc} we conclude that $(u_\eps)_\eps:=(v_\eps\ast E_\eps)_\eps$ is $(B^{{\rm{loc}}}_{p,k\wt{P_\eps}})_\eps$-moderate and defines a solution $u\in\Lb(\Gc(\R^n),\wt{\C})$ of the equation $P(D)u=v$.  
\end{proof}

\bibliographystyle{abbrv}

\begin{thebibliography}{1}


\bibitem{ChaPi:82}
J.~Chazarain and A.~Pirou.
\newblock {\em Introduction to the theory of linear partial differential
  equations}.
\newblock Number~14 in Studies in Mathematics and its Applications. North
  Holland Publishing and Co., Amsterdam, 1982.



  
  





\bibitem{FJ:98}
G.~Friedlander and M.~Joshi.
\newblock {\em Introduction to the theory of distributions}.
\newblock Cambridge University Press, New York, second edition, 1998.
 
\bibitem{Garetto:04}
C.~Garetto.
\newblock Pseudo-differential operators in algebras of generalized functions
  and global hypoellipticity.
\newblock {\em Acta Appl. Math.}, 80(2):123--174, 2004.
 
\bibitem{Garetto:05a}
C.~Garetto.
\newblock Topological structures in {C}olombeau algebras: topological
  $\wt{\C}$-modules and duality theory.
\newblock {\em Acta. Appl. Math.}, 88(1):81--123, 2005.

\bibitem{Garetto:05b}
C.~Garetto.
\newblock Topological structures in {C}olombeau algebras: investigation of the
  duals of ${\Gc(\Om)}$, ${\G(\Om)}$ and ${\GS(\R^n)}$.
\newblock {\em Monatsh. Math.}, 146(3):203--226, 2005.

\bibitem{Garetto:04th}
C.~Garetto.
\newblock Pseudodifferential operators with generalized symbols and regularity theory.
\newblock {\em Ph.D. thesis, University of Torino}, 2004.

\bibitem{Garetto:06a}
C.~Garetto.
\newblock Microlocal analysis in the dual of a {C}olombeau algebra: generalized
  wave front sets and noncharacteristic regularity.
\newblock {\em New York J. Math.}, 12:275--318, 2006.

\bibitem{Garetto:06b}
C.~Garetto.
\newblock Closed graph and open mapping theorems for topological
  $\wt{\C}$-modules and applications.
\newblock {\em arXiv:math. FA/0608087(v2)}, 2006.

\bibitem{GHO:06}
C.~Garetto, G.~H\"ormann, and M.~Oberguggenberger.
\newblock Generalized oscillatory integrals and Fourier integral operators.
\newblock {\em arXiv:math. AP/0607706}, 2006.

\bibitem{GGO:03}
C.~Garetto, T.~Gramchev, and M.~Oberguggenberger.
\newblock Pseudodifferential operators with generalized symbols and regularity
  theory.
\newblock {\em Electron. J. Diff. Eqns.}, 2005(2005)(116):1--43, 2003.

\bibitem{GH:05}
C.~Garetto and G.~H\"{o}rmann.
\newblock Microlocal analysis of generalized functions: pseudodifferential
  techniques and propagation of singularities.
\newblock {\em Proc. Edinburgh. Math. Soc.}, 48(3):603--629, 2005.

\bibitem{GH:05b}
C.~Garetto and G.~H\"ormann.
\newblock Duality theory and pseudodifferential techniques for Colombeau
  algebras: generalized kernels and microlocal analysis.
\newblock {\em Bull. Cl. Sci. Math. Nat. Sci. Math.}, 31:115--136, 2006.
 


\bibitem{HS:74}
M.~W. Hirsch and S.~Smale.
\newblock {\em Differential Equations, Dynamical Systems and Linear Algebra}.
\newblock Accademic Press, New York, 1974.

\bibitem{Hoermander:63}
L.~H{\"o}rmander.
\newblock {\em Linear Partial Differential Operators}.
\newblock Springer-Verlag, Berlin, 1963. 
 

\bibitem{Hoermander:V1}
L.~H{\"o}rmander.
\newblock {\em The Analysis of Linear Partial Differential Operators},
 volume~I.
\newblock Springer-Verlag, second edition, 1990.

\bibitem{Hoermander:V2}
L.~H{\"o}rmander.
\newblock {\em The Analysis of Linear Partial Differential Operators},
  volume~II.
\newblock Springer-Verlag, 1983.





\bibitem{HO:03}
G.~H{\"o}rmann and M.~Oberguggenberger.
\newblock Elliptic regularity and solvability for partial differential
  equations with {C}olombeau coefficients.
\newblock {\em Electron. J. Diff. Eqns.}, 2004(14):1--30, 2004.

\bibitem{HOP:05}
G.~H\"{o}rmann, M.~Oberguggenberger, and S.~Pilipovic.
\newblock Microlocal hypoellipticity of linear partial differential operators
  with generalized functions as coefficients.
\newblock {\em Trans. Amer. Math. Soc.}, 358(8):3363--3383, 2006.  












\bibitem{NP:95}
M.~Nedeljkov and S.~Pilipovi{\'{c}}.
\newblock Paley-Wiener type theorems for {C}olombeau's generalized functions.
\newblock {\em J. Math. Anal. Appl.}, 195:108--122, 1995.

\bibitem{NP:98}
M.~Nedeljkov and S.~Pilipovi{\'{c}}.
\newblock Hypoelliptic differential operators with generalized constant
  coefficients.
\newblock {\em Proc. Edinb. Math. Soc.}, 41:47--60, 1998.

\bibitem{NPS:98}
M.~Nedeljkov, S.~Pilipovi{\'{c}}, and D.~Scarpal{\'{e}}zos.
\newblock {\em The Linear Theory of Colombeau Generalized Functions}.
\newblock Pitman Research Notes in Mathematics 385. Longman Scientific {\&}
  Technical, 1998.
  


\bibitem{O:92}
M.~Oberguggenberger.
\newblock {\em Multiplication of Distributions and Applications to Partial
  Differential Equations}.
\newblock Pitman Research Notes in Mathematics 259. Longman Scientific {\&}
  Technical, 1992.
  
  






\bibitem{Scarpalezos:92}
D.~Scarpal\'{e}zos.
\newblock Topologies dans les espaces de nouvelles fonctions
  g\'{e}n\'{e}ralis\'{e}es de {C}olombeau. ${\widetilde{\C}}$-modules
  topologiques.
\newblock Universit\'{e} Paris 7, 1992.

\bibitem{Scarpalezos:98}
D.~Scarpal\'{e}zos.
\newblock Some remarks on functoriality of {C}olombeau's construction;
  topological and microlocal aspects and applications.
\newblock {\em Integral Transform. Spec. Funct.}, 6(1-4):295--307, 1998.

\bibitem{Scarpalezos:00}
D.~Scarpal\'{e}zos.
\newblock Colombeau's generalized functions: topological structures; microlocal
  properties. a simplified point of view.
\newblock {\em I. Bull. Cl. Sci. Math. Nat. Sci. Math.}, 25:89--114, 2000.

\bibitem{Soraggi:96}
R.~L. Soraggi.
\newblock Fourier analysis on {C}olombeau's algebra of generalized functions.
\newblock {\em J. Anal. Math}, 69:201--227, 1996.


\end{thebibliography}
\newcommand{\SortNoop}[1]{}

\end{document}